\documentclass[reqno]{amsart}
\usepackage{amssymb}
\usepackage{amsthm}
\usepackage{amsmath}
\usepackage{enumitem}
\usepackage{xcolor}
\usepackage{booktabs}
\usepackage[foot]{amsaddr}

\usepackage[hyphens]{url}
\usepackage{hyperref}
\usepackage[hyphenbreaks]{breakurl}

\usepackage{orcidlink}
\usepackage{subcaption}

\hypersetup{
    colorlinks=true,
    linkcolor=babyblue,
    filecolor=babyblue,      
    urlcolor=babyblue,
    citecolor=babyblue
    }

\definecolor{babyblue}{rgb}{0.1, 0.6, 0.75}

\definecolor{orcidlogocol}{rgb}{0.1, 0.6, 0.75}

\newtheorem{theorem}{Theorem}[section]
\newtheorem{corollary}[theorem]{Corollary}
\newtheorem{proposition}[theorem]{Proposition}
\newtheorem{lemma}[theorem]{Lemma}

\theoremstyle{definition}

\def\Aut{{\rm Aut}}

\def\Sym{{\rm Sym}}

\def\ker{{\rm ker}\,}
\def\ord{{\rm ord}}

\def\Red{{\rm Red}}

\def\Comp{{\rm Comp}}

\newcommand{\NN}{\mathbb N}
\newcommand{\ZZ}{\mathbb Z}

\def\CC{{C\nolinebreak[4]\hspace{-.05em}\raisebox{.4ex}{\tiny\bf ++}}}

\begin{document}
\title{Recursive characterisation of skew morphisms of finite cyclic groups}
%\date{\today}
\author[Martin Bachrat{\' y}]{Martin Bachrat{\' y}\,$^{\orcidlink{0000-0002-4300-7507},1,\ast}$}
\address{$^1$Faculty of Civil Engineering, Slovak University of Technology, Bratislava 81005, Slovakia}
\address{$^\ast$Corresponding author}
\email[A1]{martin.bachraty@stuba.sk}
\thanks{The first author acknowledges funding from the EU NextGenerationEU through the Recovery and Resilience Plan for Slovakia under the project No. 09I03-03-V04-00272.}
\author[Michal Hagara]{Michal Hagara}
%\address{$^2$Bratislava, Slovakia}
\email[A2]{michal.hagara@gmail.com}

\maketitle

\begin{abstract}
A skew morphism of a finite group $G$ is an element $\varphi$ of $\Sym(G)$ preserving the identity element of $G$ and having the property that for each $a\in G$ there exists a non-negative integer $i_a$ such that $\varphi(ab)=\varphi(a)\varphi^{i_a}(b)$ for all $b\in G$. In this paper we show that if a skew morphism $\varphi$ of $\ZZ_n$ is not an automorphism of $\ZZ_n$, then it is uniquely determined by a triple $(h,\alpha,\beta)$ where $h$ is an element of $\ZZ_n$, $\alpha$ is a skew morphism of $\ZZ_a$ where $a<n$, and $\beta$ is a skew morphism of $\ZZ_b$ where either $b<n$, or $b=n$ and $|\langle \beta\rangle| <|\langle \varphi\rangle|$. Conversely, we also list necessary and sufficient conditions for a triple $(h,\alpha,\beta)$ to define a skew morphism of a given cyclic group. In particular, this gives a recursive characterisation of skew morphisms for all finite cyclic groups. We use this characterisation to prove new theorems about skew morphisms of cyclic groups and to generate a census of all skew morphisms for cyclic groups of order up to $2000$. \\[0.4em]
\end{abstract}

%%%%%%%%%%%
\section{Introduction}\label{sec:prel}
%%%%%%%%%%%

A permutation $\varphi$ of a finite group $G$ is called a \emph{skew morphism} of $G$ if it fixes the identity of $G$ and for each $a\in G$ there exist a non-negative integer $i_a$ such that $\varphi(ab)=\varphi(a)\varphi^{i_a}(b)$ for all $b\in G$. The \emph{order} of $\varphi$ is defined as the order of the cyclic group $\langle \varphi \rangle$, and denoted by $\ord(\varphi)$. It is straightforward to check that if $\varphi$ is non-trivial, then for each $a\in G$ there is a unique choice for $i_a$ such that $i_a\in \{1, \dots, \ord(\varphi)-1\}$. The function $\pi_{\varphi}$ that maps each $a$ to this integer $i_a$ is called the \emph{power function} of $\varphi$. In the case when $\varphi$ is the trivial permutation of $G$, we define $\pi_{\varphi}(a)=1$ for all $a\in G$. Note that every automorphism of $G$ is also a skew morphism of $G$. We will refer to skew morphisms that are not automorphisms as \emph{proper} skew morphisms.    

The concept of a skew morphism was first introduced more than two decades ago in~\cite{JajcaySiran}. At that time the primary interest in skew morphisms came from their connection to regular Cayley maps (see~\cite{ConderJajcayTucker2007b, ConderJajcayTucker2007, FengJajcayWang} for example), but soon enough the study of skew morphisms became also an independent topic in group theory (see~\cite{BachratyConderVerret,ConderJajcayTucker2016,HuKovacsKwon} for example). One of the most prominent problems in this field is to classify skew morphisms for all finite cyclic groups. Despite substantial progress over the last decade this problem still remains open. Interestingly, skew morphisms have been already classified for other important families of finite groups such as dihedral groups or finite simple groups.

Before we proceed let us define a strict partial order on the set of all skew morphisms of finite cyclic groups. For skew morphisms $\varphi$ of $\ZZ_n$ and $\rho$ of $\ZZ_m$ we say that $\varphi$ is \emph{smaller} than $\rho$ if $n<m$, or $n=m$ and $\ord(\varphi)<\ord(\rho)$. In this paper we provide a recursive characterisation of skew morphisms for all finite cyclic groups by showing that each skew morphism $\varphi$ of $\ZZ_n$ is uniquely determined by an element $h$ of $\ZZ_n$ and a pair of skew morphisms $\alpha$ and $\beta$ (of some finite cyclic groups) such that both $\alpha$ and $\beta$ are smaller than $\varphi$. We also explain how to quickly identify all triples $(h,\alpha,\beta)$ that actually determine some skew morphism of $\ZZ_n$. By ``quickly'' here we mean that the time complexity of the algorithm that checks a given triple and constructs the corresponding skew morphism is linear in $n$. Using this method we can easily find all skew morphisms of a given finite cyclic group, provided that we already have a complete list of skew morphisms for all cyclic groups of smaller orders.

In Section~\ref{sec:prel} we give some necessary background on skew morphisms of finite groups. Then in Section~\ref{sec:red} we discuss different possible ways of reducing a skew morphism of a finite cyclic group to a smaller skew morphism, and we formulate our main theorem. The proof of this theorem presented in Section~\ref{sec:main} is followed by Section~\ref{sec:recs} in which we explain in detail the recursive method for constructing all skew morphisms of cyclic groups up to any given finite order. In Section~\ref{sec:apps} we show theoretical applications of our prior observations, and this is followed by concluding remarks in Section~\ref{sec:remarks}.

%%%%%%%%%%%
\section{Preliminaries}\label{sec:prel}
%%%%%%%%%%%

Let $\varphi$ be a skew morphism of a finite group $G$ with power function $\pi_{\varphi}$, and let $\ker\varphi = \{a\in G \mid \pi_{\varphi}(a) = 1 \}$. The set $\ker\varphi$ is called the \emph{kernel} of $\varphi$, and it can be easily shown that it is a subgroup of $G$. In the case when $G$ is abelian, the kernel is preserved by $\varphi$ set-wise. Any two elements $a,b\in G$ are in the same right coset of $\ker\varphi$ in $G$ if and only if $\pi_{\varphi}(a)=\pi_{\varphi}(b)$; see~\cite{ConderJajcayTucker2007} for example. If $G=\ZZ_n$, we have $\pi_{\varphi}(a)=\pi_{\varphi}(b)$ if and only if $a\equiv b \pmod{n/|\ker\varphi|}$. We will repeatedly use the fact (proved in~\cite{ConderJajcayTucker2016}) that every skew morphism of a non-trivial group has a non-trivial kernel. If $p$ is a prime, then the only non-trivial subgroup of $\ZZ_p$ is $\ZZ_p$ itself, and so $\ZZ_p$ admits no proper skew morphisms. Next, note that for each $a\in G$ we can view $\pi_{\varphi}$ as an element of $\ZZ_{\ord(\varphi)}$. In what follows we will denote the element $\sum_{0\leq i \leq x-1}\pi_{\varphi}(\varphi^i(y))$ of $\ZZ_{\ord(\varphi)}$ by $\sigma_{\varphi}(x,y)$. Also note that if $\varphi(H)=H$ for some subgroup $H$ of $G$, then from the definition of a skew morphism it follows that the restriction of $\varphi$ to $H$ is a skew morphism of $H$. This skew morphism will be denoted by $\varphi\vert_H$. The following facts about skew morphisms will be useful.

\begin{lemma}[\cite{ConderJajcayTucker2016}]\label{lem:powprod}
Let $\varphi$ be a skew morphism of a finite group $G$ with power function $\pi_{\varphi}$, and let $e$ be a positive integer. Then $\varphi^e(ab)=\varphi^e(a)\varphi^{\sigma_{\varphi}(e,a)}(b)$ and $\pi_{\varphi}(ab)=\sigma_{\varphi}(\pi_{\varphi}(a),b)$ for every $a,b\in G$.
\end{lemma}

\begin{theorem}[\cite{ConderJajcayTucker2016}]
\label{thm:orderofskew}
The order of a skew morphism of a non-trivial finite group $G$ is less than the order of $G$.
\end{theorem} 

\begin{theorem}[\cite{KovacsNedela2011}]
\label{thm:orderofskewcyclic}
If $\varphi$ is a skew morphism of a group $\ZZ_n$, then $\ord(\varphi)$ divides $n\phi(n)$. Moreover, if $\gcd(\ord(\varphi),n)=1$, then $\varphi$ is an automorphism of $\ZZ_n$.
\end{theorem}

\begin{theorem}[\cite{KovacsNedela2011}]
\label{thm:orderoforbit}
Let $\varphi$ be a skew morphism of a finite group $G$, and let $T$ be an orbit of $\langle \varphi \rangle$. If $T$ generates $G$, then $\ord(\varphi)=|T|$. In particular, if $G=\ZZ_n$, then $\ord(\varphi)$ is equal to the smallest positive integer $i$ satisfying $\varphi^i(1)=1$.
\end{theorem} 

\begin{lemma}[\cite{BachratyJajcay2016}]
\label{lem:power}
Let $\varphi$ be a skew morphism of a finite group $G$, and let $i$ be a positive integer. Then $\varphi^i$ is a skew morphism of $G$ if and only if for every $a\in G$ there exists some $x_a\in \ZZ_{\ord(\varphi)}$ such that $\sigma_{\varphi}(i,a)\equiv ix_a \pmod{\ord(\varphi)}$. If $\varphi^i$ is a skew morphism of $G$, then $\pi_{\varphi^i}(a)$ is the smallest positive solution $x_a$ of the above equation.
\end{lemma}

\begin{lemma}[\cite{Bachraty}]\label{lem:skewfromorbit}
Let $\varphi$ be a skew morphism of $\ZZ_n$ with power function $\pi_{\varphi}$. Then:
\begin{align*}
\varphi(i) = \varphi(1) + \varphi^{\pi_{\varphi}(1)}(1) + \varphi^{\pi_{\varphi}(2)}(1) +  \dots +  \varphi^{\pi_{\varphi}(i-1)}(1) \hbox{ for all } i\in \NN .
\end{align*}
\end{lemma}

\begin{lemma}[\cite{JajcaySiran}]\label{lem:star}
Let $\varphi$ be a skew morphism of $\ZZ_n$ with kernel $\ker\varphi$, and let $k=n/|\ker\varphi|$. Then the mapping $\varphi^*\!: \ZZ_n/(\ker\varphi) \to \ZZ_n/(\ker\varphi)$ given by $\varphi^*(a) = \varphi(a) \bmod{k}$ is a well-defined skew morphism of $\ZZ_n/(\ker\varphi)$.
\end{lemma}

Let $\varphi$ be a skew morphism of $\ZZ_n$, and let $m$ be a divisor of $n$ such that for any pair $a,b\in \ZZ_n$ we have $a\equiv b \pmod{m}$ if and only if $\varphi(a)\equiv \varphi(b) \pmod{m}$. Then the skew morphism $\varphi$ \emph{taken modulo} $m$ is the permutation of $\ZZ_m$ defined by $x\mapsto \varphi(x) \bmod{m}$. Note that the skew morphism $\varphi^*$ described in Lemma~\ref{lem:star} can be defined as $\varphi$ taken modulo $n/|\ker\varphi|$. 

\begin{lemma}[\cite{KovacsNedela2011}]\label{lem:derived}
If $\varphi$ is a skew morphism of $\ZZ_n$, then the mapping $\varphi' \!: \ZZ_{\ord(\varphi)} \to \ZZ_{\ord(\varphi)}$ given by $\varphi'(a) = \sigma_{\varphi}(a,1)$ is a well-defined skew morphism of $\ZZ_{\ord(\varphi)}$.
\end{lemma}

The skew morphism $\varphi'$ of $\ZZ_{\ord(\varphi)}$ described in Lemma~\ref{lem:derived} is called the skew morphism \emph{derived} from $\varphi$. The concept of the derived skew morphism (sometimes also called the quotient of a skew morphism) was first introduced in~\cite{KovacsNedela2011}, and it was later studied in more detail in~\cite{Bachraty, BachratyHagara, FengHu, HuNedela}. Derived skew morphisms will play a key role in later sections, here we list some of their useful properties.

\begin{proposition}[\cite{BachratyHagara}]\label{prop:quo}
Let $\varphi$ be a non-trivial skew morphism of $\ZZ_n$, and let $\varphi'$ be the skew morphism derived from $\varphi$. Then:
\begin{enumerate}[label={\rm (\alph*)},ref=\ref{prop:quo}(\alph*)]
\item\label{prop:quo:a} $(\varphi')^b(a) = \sigma_{\varphi}(a,b)$; 
\item\label{prop:quo:b} $\pi_{\varphi}(i) = (\varphi')^i(1)$, so in particular $\ord(\varphi')=n/|\ker\varphi|$; 
\item\label{prop:quo:c} $\varphi^i(1) \equiv \pi_{\varphi'}(i) \pmod{n/|\ker\varphi|}$;
\item\label{prop:quo:d} $\ord(\varphi)= \ord(\varphi^*)|\ker\varphi'|$;
\item\label{prop:quo:e} the skew morphism derived from $\varphi'$ is equal to $\varphi^*$;
\item\label{prop:quo:f} $\varphi$ is an automorphism if and only if $\varphi'$ is the identity mapping; and
\item\label{prop:quo:g} $\varphi$ preserves the cosets of $\ker\varphi$ in $\ZZ_n$ if and only if $\varphi'$ is an automorphism.
\end{enumerate}
\end{proposition}

\begin{theorem}[\cite{BachratyHagara}]\label{thm:power}
Let $\varphi$ be a skew morphism of $\ZZ_n$, and let $e$ be a positive integer. Then $\varphi^e$ is a skew morphism of $\ZZ_n$ if and only if the subgroup of $\ZZ_{\ord(\varphi)}$ generated by $e$ is preserved by $\varphi'$ set-wise. 
\end{theorem}

An important consequence of Theorem~\ref{thm:power} (proved in~\cite{BachratyHagara}) is the fact that if $p$ is the smallest prime divisor of the order of a non-trivial skew morphism $\varphi$ of $\ZZ_n$, then $\varphi^p$ is a skew morphism of $\ZZ_n$. 

%Next, we say that two skew morphism $\varphi_1$ and $\varphi_2$ of $\ZZ_n$ are weakly equivalent if and only if $\varphi_2=\delta^{-1}(\varphi_1)^e\delta$ for some $\delta\in \Aut(\ZZ_n)$ and some positive integer $e$ relatively prime to $\ord(\varphi_1)$. The following observations about the weak equivalence will be helpful.

%\begin{lemma}[\cite{BachratyHagara}]\label{lem:weak}
%The weak equivalence is an equivalence relation. If two skew morphisms $\varphi_1$ and $\varphi_2$ of a finite group $G$ are weakly equivalent, then:
%\begin{enumerate}[label={\rm (\alph*)},ref=\ref{lem:weak}(\alph*)]
%\item\label{lem:weak:a} $\ord(\varphi_1)=\ord(\varphi_2)$; 
%\item\label{lem:weak:b} $\ker\varphi_1=\ker\varphi_2$; 
%\item\label{lem:weak:c} if either $\varphi_1^{\,e}$ or $\varphi_2^{\,e}$ is a skew morphism of $G$, then so is the other and the two skew morphisms are weakly equivalent; and
%\item\label{lem:weak:d} if $G=\ZZ_n$, then $(\varphi_1)'$ and $(\varphi_2)'$ are weakly equivalent. 
%\end{enumerate}
%\end{lemma}    

%%%%%%%%%%%
\section{Reducibility of skew morphisms}\label{sec:red}
%%%%%%%%%%%

Let $\varphi$ be a skew morphism of $\ZZ_n$. By a \emph{reduction} of $\varphi$ we will understand any unary operation that maps $\varphi$ to a skew morphism $\rho$ of $\ZZ_r$ such that either $r<n$, or $r=n$ and $\ord(\rho)<\ord(\varphi)$. We have already seen some examples of reductions in Section~\ref{sec:prel}. Namely, unless $n=1$ both $\varphi'$ and $\varphi^*$ are reductions of $\varphi$. Further, we know that $\ker\varphi$ is preserved by $\varphi$, and so if $|\ker\varphi|<n$ (or, in other words, if $\varphi$ is proper), then $\varphi\vert_{\ker\varphi}$ is also a reduction of $\varphi$. To provide another example note that if the order of $\varphi$ is non-trivial, then $\ord(\varphi)$ has the smallest prime divisor, say $p$, and $\varphi^p$ is a reduction of $\varphi$.

It can be easily checked that all of the above reductions are feasible for all proper skew morphisms of cyclic groups. In contrast, some reductions might not be always viable. An example of such reduction can be found in~\cite{KovacsNedela2011} where it was shown that if $n=ab$ for some positive integers $a$ and $b$ such that $\gcd(a,b)=\gcd(\phi(a),b)=\gcd(a,\phi(b))=1$, then $\varphi$ taken modulo $a$ and modulo $b$ restricts to skew morphisms $\varphi_a$ and $\varphi_b$ of $\ZZ_a$ and $\ZZ_b$, respectively. If both $a$ and $b$ are strictly smaller than $n$, then both $\varphi_a$ and $\varphi_b$ are reductions of $\varphi$. Interestingly, these two reductions can be used to reconstruct $\varphi$. Namely, since $\gcd(a,b)=1$, it follows that $\ZZ_n = \ZZ_a \times \ZZ_b$, and hence $\varphi = \varphi_a \times \varphi_b$. The main issue of this method of constructing skew morphisms is that there are infinitely many values $n$ for which the only suitable decomposition is $n\cdot 1$ (or $1 \cdot n$), in which case $\varphi_a$ (or $\varphi_b$) is not a reduction.

We already noted that if $\varphi$ is proper, then $\varphi'$, $\varphi^p$ and $\varphi\vert_{\ker\varphi}$ are all well-defined reductions of $\varphi$.  (We can also add a fourth reduction $\varphi^*$, but this is redundant since by Proposition~\ref{prop:quo:e} we know that $\varphi^*=(\varphi')'$.) In contrast to the previous case (where we reduced $\varphi$ to $\varphi_a$ and $\varphi_b$), the original skew morphism $\varphi$ is not uniquely determined by these three reductions. The smallest counterexample (in terms of $n$) is the pair of 
%(weakly equivalent) 
skew morphisms $(1,3,5)$ and $(1,5,3)$ of $\ZZ_6$ with kernels of order $3$. In both cases the derived skew morphisms are equal to the automorphism $(1,2)$ of $\ZZ_3$, restrictions to the kernel are equal to the trivial permutation of $\ZZ_3$, and $p$\,th powers are equal to the trivial permutation of $\ZZ_6$. %The smallest counterexample where the two skew morphisms (with the three reductions being equal) are not weakly equivalent exists for $n=72$ {\color{red}IS THIS TRUE?}.

The aim of this paper is to find a collection of reductions feasible for all proper skew morphisms (of cyclic groups) that will always uniquely determine the original skew morphism. This will give a recursive characterisation of all skew morphisms for all finite cyclic groups. The main theorem of this paper is the following:

\begin{theorem}\label{thm:main}
Let $\varphi$ be a proper skew morphism of $\ZZ_n$. Then there exists an element $h\in \ZZ_n$ and a pair of reductions $\alpha$ and $\beta$ of $\varphi$ that uniquely determine $\varphi$. 
\end{theorem} 

The proof of Theorem~\ref{thm:main} is presented in Section~\ref{sec:main}. Then in Section~\ref{sec:recs} we explain how the recursive characterisation can be used to find all skew morphisms for cyclic groups up to any order. Finally in Section~\ref{sec:apps} we look at various applications of our findings.

%In what follows we will often refer to two general types of reductions. If $m$ is a proper divisor of $n$ such that $m\geq 2$ and $\varphi(\langle m \rangle) = \langle m \rangle$, then clearly $\varphi\vert_{\langle m \rangle}$ is a reduction of $\varphi$, and we say that $\varphi$ is \emph{$m$-reducible to a subgroup}. If $m$ is a proper divisor of $\ord(\varphi)$ such that $m\geq 2$ and $\varphi^m$ is a skew morphism of $\ZZ_n$, then $\varphi^m$ is a reduction of $\varphi$, and we say that $\varphi$ is \emph{$m$-reducible to a power}. 

%%%%%%%%%%%
\section{Proof of Theorem~\ref{thm:main}}\label{sec:main}
%%%%%%%%%%%

Let $\varphi$ be a skew morphism of a non-trivial cyclic group. Recall that $\varphi'$ is a skew morphism of $\ZZ_{\ord(\varphi)}$, and so by Theorem~\ref{thm:orderofskew} it follows that by taking the derived skew morphism of $\varphi$ we decrease the order of the underlying group. Hence by repeating the operation of taking the derived skew morphism we will eventually reduce $\varphi$ to the identity permutation of the trivial group. Since we start with a skew morphism of a non-trivial group, this reduction process must have the penultimate step, and so $\varphi$ must be eventually reduced to a skew morphism $\rho$ of some non-trivial cyclic group such that $\rho'$ is the identity permutation of the trivial group. Hence $\rho'$ is a skew morphism of $\ZZ_1$, which implies $\ord(\rho)=1$, and it follows that $\rho$ is the trivial permutation of some non-trivial cyclic group. Let $m$ denote the order of this non-trivial cyclic group, and let $c$ denote the number of steps needed for reducing $\varphi$ to $\rho$ (by taking derived skew morphisms). We say that $\varphi$ has \emph{complexity} $c$ and \emph{auto-order} $m$. For every non-negative integer $i$ we also let $\varphi^{(i)}$ denote $\varphi$ reduced $i$ times; in particular, $\varphi^{(0)}=\varphi$, $\varphi^{(1)}=\varphi'$, and $\varphi^{(c)}$ is the identity permutation of $\ZZ_m$. Also note that by Proposition~\ref{prop:quo:e} we have $\varphi^{(2)}=\varphi^*$. In the following lemma we list some useful properties of complexities and auto-orders. 

\begin{lemma}\label{lem:comp}
Let $n\geq 2$, and let $\varphi$ be a skew morphism of $\ZZ_n$ with complexity $c$ and auto-order $m$. Then:
\begin{enumerate}[label={\rm (\alph*)},ref=\ref{lem:comp}(\alph*)]
\item\label{lem:comp:a} $c\geq 0$ and $m\geq 2$;
\item\label{lem:comp:b} $\varphi$ is the trivial permutation if and only if $c=0$; 
\item\label{lem:comp:c} $\varphi$ is a non-trivial automorphism if and only if $c=1$;
\item\label{lem:comp:d} $\varphi$ is a proper skew morphism if and only if $c\geq 2$;
\item\label{lem:comp:e} if $c\geq 1$, then $\varphi'$ has complexity $c-1$ and auto-order $m$;
\item\label{lem:comp:f} if $c\geq 1$, then $\varphi^{(c-1)}$ is an automorphism of order $m$. 
%\item\label{lem:comp:g} if $c\geq 1$, then $\varphi^{(c)}$ is the identity permutation of $\ZZ_{\ord(\varphi^{(c-1)})}$.
%\item\label{lem:comp:h} weakly equivalent skew morphisms have the same complexity and auto-order.
\end{enumerate}
\end{lemma} 
\begin{proof}
The proof follow trivially from the definition of the complexity and Proposition~\ref{prop:quo}. %To prove (h) note that by Lemma~\ref{lem:weak:d} we know that if two skew morphisms are weakly equivalent then so are their derived skew morphisms. It follows that if some skew morphism $\rho$ of $\ZZ_n$ is weakly equivalent to $\varphi$, then $\rho^{(c)}$ is weakly equivalent to $\varphi^{(c)}$, which is the identity permutation of $\ZZ_m$. Hence by Lemma~\ref{lem:weak:b} we have $\ord(\rho^{(c)})=\ord(\varphi^{(c)})=1$, and so $\rho^{(c)}$ is the identity permutation of $\ZZ_m$ too. Therefore $\rho$ has complexity $c$ and auto-order $m$. 
\end{proof}   

Next we prove a useful observation about auto-orders of skew morphisms with even complexity.

\begin{lemma}\label{lem:main}
Let $\varphi$ be a skew morphism of $\ZZ_n$ with even complexity and auto-order $m$. Then $m$ divides $n$, and $\varphi$ taken modulo $m$ is the trivial permutation of $\ZZ_m$.
\end{lemma}
\begin{proof}
For each $i\in \{0,\dots, c\}$ let $n_i$ denote the positive integer such that $\varphi^{(i)}$ is a skew morphism of $\ZZ_{n_i}$. Also let $c$ be the complexity of $\varphi$, hence $c$ is even and $\varphi^{(c)}$ is the identity permutation of $\ZZ_m$. Recall that for every non-negative integer $i$ we have $\varphi^{(i+2)} = (\varphi^{(i)})^*$, and so by Lemma~\ref{lem:star} we have $n_{i+2}=n_i/|\ker (\varphi^{(i)})|$, and hence $n_{i+2}$ divides $n_i$. In particular, it follows that $n_c$ divides $n_{c-2}$, which divides $n_{c-4}$, which divides $n_{c-6}$, and so on. Then, since $c$ is even, we deduce that $n_c$ divides $n_0$, and the first assertion follows from the fact that $n_c=m$ and $n_0=n$.  

To prove the second part of the assertion note that $\varphi^{(i+2)}$ (which is equal to $(\varphi^{(i)})^*$) can be viewed as $\varphi^{(i)}$ taken modulo $n_{i+2}$. Hence, if $m$ divides $n_{i+2}$ and $\varphi^{(i+2)}$ taken modulo $m$ restricts to the identity of $\ZZ_m$, then the same is true also for $\varphi^{(i)}$. (Here we also use the observation that since $m$ divides $n_{i+2}$, it also divides $n_i$.) Since $\varphi^{(c)}$ is the trivial permutation of $\ZZ_m$ and $c$ is even, the proof of the second assertion now follows easily by induction.
\end{proof}

We are now ready to present the proof of Theorem~\ref{thm:main}. It splits into two cases depending on the parity of the complexity of a skew morphism.

\begin{theorem}\label{thm:even}
Let $\varphi$ be a skew morphism of $\ZZ_n$ with even complexity $c$ and auto-order $m$. Then $\varphi$ is uniquely determined by the triple $(\varphi(1),\varphi',\varphi\vert_{\langle m\rangle})$.
\end{theorem}
\begin{proof}
Note that by Lemma~\ref{lem:main} we know that $m$ divides $n$ and $\varphi$ taken modulo $m$ is the trivial permutation of $\ZZ_m$. It follows that $\langle m \rangle$ is a subgroup of $\ZZ_n$, and since elements of this subgroup are exactly the elements of $\ZZ_n$ equal to $0$ modulo $m$, we have $\varphi(\langle m\rangle)=\langle m \rangle$. It follows that $\varphi\vert_{\langle m\rangle}$ is equal to some skew morphism $\beta$ of $\ZZ_{|\langle m \rangle|}$ and $\varphi(im)=\beta(i)m$ for each $i\in \ZZ_{n/m}$. Another consequence of the fact that $\varphi$ taken modulo $m$ is the trivial permutation of $\ZZ_m$ is that $\varphi(1) \equiv 1 \pmod{m}$, and so there exists some integer $f$ such that $\varphi(1)=1+fm$. It follows that $\varphi^2(1)=\varphi(\varphi(1))=\varphi(1+fm)=\varphi(1)+\varphi^{\pi_{\varphi}(1)}(fm)=\varphi(1)+\beta^{\pi_{\varphi}(1)}(f)m$, and by induction it can be easily checked that $\varphi^e(1)=\varphi^{e-1}(1)+\beta^{\pi_{\varphi}(1)+\pi_{\varphi}(\varphi(1))+\dots+\pi_{\varphi}(\varphi^{e-2}(1))}(f)m$ for all $e\geq 2$. Note that by definitions of $\sigma_{\varphi}$ and $\varphi'$ we have $\pi_{\varphi}(1)+\pi_{\varphi}(\varphi(1))+\dots+\pi_{\varphi}(\varphi^{e-2}(1)) \equiv \sigma_{\varphi}(e-1,1) \equiv \varphi'(e-1) \pmod{\ord(\varphi)}$, and since clearly $\ord(\beta)$ divides $\ord(\varphi)$, we have $\varphi^e(1) = \varphi^{e-1}(1) + \beta^{\varphi'(e-1)}(f)m$. Hence, if we let $\alpha = \varphi'$, it follows that 
\begin{equation*}
\varphi^e(1)=1+fm+\beta^{\alpha(1)}(f)m+\dots + \beta^{\alpha(e-1)}(f)m \text{ for all } e\geq 2. 
\end{equation*}
Next, by Lemma~\ref{lem:comp:e} we know that $m$ is equal to the auto-order of $\alpha$, and we also have $f=(\varphi(1)-1)/m$. It follows that the orbit of $\langle \varphi \rangle$ containing $1$ is uniquely determined by $\varphi(1)$, $\alpha$ and $\beta$. The rest follows by Lemma~\ref{lem:skewfromorbit} and Proposition~\ref{prop:quo:b} which implies that $\pi_{\varphi}(i) = \alpha^i(1)$ for all $i\in \ZZ_n$.   
\end{proof}

\begin{theorem}\label{thm:odd}
Let $\varphi$ be a skew morphism of $\ZZ_n$ with odd complexity $c$ and auto-order $m$. Then $\varphi$ is uniquely determined by the triple $(\varphi(1),\varphi',\varphi^m)$. 
\end{theorem}
\begin{proof}
By Lemma~\ref{lem:comp:e} we know that $\varphi'$ has (even) complexity $c-1$ and auto-order $m$, and consequently by Lemma~\ref{lem:main} we deduce that the subgroup of $\ZZ_{\ord(\varphi)}$ generated by $m$ is preserved by $\varphi'$ set-wise. Hence by Theorem~\ref{thm:power} we know that $\varphi^m$ is a skew morphism of $\ZZ_n$. Furthermore, we know (again by Lemma~\ref{lem:main}) that $\varphi'$ taken modulo $m$ is the trivial permutation of $\ZZ_m$, and combining this with Proposition~\ref{prop:quo:b} we deduce that $\pi_{\varphi}(i)\equiv (\varphi')^i(1) \equiv 1 \pmod{m}$ for all $i\in \ZZ_n$. It follows that $\pi_{\varphi}(i)-1$ is always a multiple of $m$ or, equivalently, that $(\pi_{\varphi}(i)-1)/m$ is an integer. Let $h=\varphi(1)$, $\alpha=\varphi'$ and $\beta=\varphi^m$, and recall that $\pi_{\varphi}(i)=\alpha^i(1)$. Also note that $\varphi^{\pi_{\varphi}(i)}(1)=\varphi^{\pi_{\varphi}(i)-1}(\varphi(1))=\beta^{(\pi_{\varphi}(i)-1)/m}(h)=\beta^{(\alpha^i(1)-1)/m}(h)$ for all $i\in \ZZ_n$, and then by Lemma~\ref{lem:skewfromorbit} we have 
\begin{equation*}
\varphi(a) = h + \beta^{(\alpha(1)-1)/m}(h) + \dots + \beta^{(\alpha^{a-1}(1)-1)/m}(h) \text{ for all } a\in \ZZ_n.
\end{equation*}
Since $m$ is equal to the auto-order of $\alpha$, it follows that $\varphi$ is uniquely determined by $h$, $\alpha$ and $\beta$.
\end{proof}

For a skew morphism $\varphi$ of $\ZZ_n$ with complexity $c$ and auto-order $m$ we let $\Red(\varphi)$ denote the triple $(h,\alpha,\beta)$ where $h=\varphi(1)$, $\alpha$ is the skew morphism derived from $\varphi$, and $\beta$ is either the restriction of $\varphi$ to $\langle m\rangle$ or the $m$\,th power of $\varphi$, depending on whether $c$ is even or odd. Recall that $\alpha$ is a reduction of $\varphi$, as argued in the previous section. By Lemma~\ref{lem:comp:a} we have $m\geq 2$, and so by Lemma~\ref{lem:main} we deduce that if $c$ is even, then $\langle m\rangle$ is a proper subgroup of $\ZZ_n$. This proves that if $c$ is even, then $\beta$ is a reduction of $\varphi$.
For $c$ odd we showed in the proof of Theorem~\ref{thm:odd} that $m$ divides $\ord(\varphi)$, which together with the fact that $m\geq 2$ implies that $\ord(\varphi^m)<\ord(\varphi)$, and hence $\beta$ is a reduction of $\varphi$ in this case as well. Finally, an immediate consequence of Theorems~\ref{thm:even} and~\ref{thm:odd} is that $\Red(\varphi)$ uniquely determines $\varphi$, which completes the proof of Theorem~\ref{thm:main}.

%%%%%%%%%%%
\section{Recursive construction of skew morphisms}\label{sec:recs}
%%%%%%%%%%%

Using Theorem~\ref{thm:main} we can find recursively all skew morphisms of cyclic groups up to any given order $N$. The method goes as follows: 
\begin{enumerate}[label={\rm (\arabic*)},ref=(\arabic*)]
\item\label{step:1} It is well-known that the order of the smallest cyclic group that admits a proper skew morphism is $6$, so we set $n=6$. 
\item\label{step:2} The smallest possible order of a proper skew morphism of $\ZZ_n$ is $3$ (this is an easy consequence of Proposition~\ref{prop:quo:f}), so we set $\ell=3$.
\item\label{step:3} For every proper skew morphism $\varphi$ of $\ZZ_n$ of order $\ell$ there exists the reduction $\Red(\varphi)$ to an element of $\ZZ_n$ and two skew morphisms $\alpha$ and $\beta$ that were both already discovered. It follows that we can list all possible candidates $(h,\alpha,\beta)$ for $\Red(\varphi)$. Note that $\alpha$ or $\beta$ might be automorphisms, but these are all known for cyclic groups.  
\item\label{step:4} For each candidate $(h,\alpha,\beta)$ we check if it defines a skew morphism of $\ZZ_n$ of order $\ell$.
\item\label{step:5} We increase the value of $\ell$ by one. If $\ell<n$, we return to Step~\ref{step:3}. If $\ell=n$, we increase the value of $n$ by one and return to Step~\ref{step:2}. This is repeated until $n$ exceeds $N$.   
\end{enumerate}
When applying this method in practice, the most crucial parts are Step~\ref{step:3} and Step~\ref{step:4}. In what follows we will show how to limit the number of possible candidates in Step~\ref{step:3}, and we will prove that the check in Step~\ref{step:4} and consequent construction of $\varphi$ can be carried out in linear time with respect to $n$. We will consider two separate cases depending on the parity of the complexity of $\varphi$.

\subsection{Skew morphisms with odd complexity}%%%%%%    
Let $\varphi$ be a skew morphism of $\ZZ_n$ with odd complexity and auto-order $m$, and let $\Red(\varphi)=(h,\alpha,\beta)$. By Lemma~\ref{lem:comp:e} we know that $\alpha$ has even complexity and auto-order $m$, and so by Lemma~\ref{lem:main} we have $\alpha^e(1)\equiv 1 \pmod{m}$ for every positive integer $e$. It follows that $\alpha^e(1)-1$ is divisible by $m$, and for each $e$ we define the integer $k_e$ by $k_e=(\alpha^e(1)-1)/m$. Then from Theorem~\ref{thm:odd} we deduce that
\begin{equation}\label{eq:odd}
\varphi(a) = h + \beta^{k_1}(h) + \dots + \beta^{k_{a-1}}(h) \text{ for all } a\in \ZZ_n\, .
\end{equation}
In the following proposition we characterise all triples $(h,\alpha,\beta)$ such that Equation~\eqref{eq:odd} defines a skew morphism of $\ZZ_n$. Since Equation~\eqref{eq:odd} is well-defined only for positive integers, we define $\varphi(0)$ as $\varphi(n)$. 
\begin{proposition}\label{prop:odd}
Let $n$ and $m$ be positive integers, let $h$ be an element of $\ZZ_n$, let $\beta$ be a skew morphism of $\ZZ_n$, and let $\alpha$ be a skew morphism of $\ZZ_{\ord(\beta)m}$ with even complexity and auto-order $m$ such that $\ord(\alpha)$ divides $n$. Also let $\varphi$ be the mapping from $\ZZ_n$ to $\ZZ_n$ given by Equation~\eqref{eq:odd}. Then $\varphi$ is a skew morphism of $\ZZ_n$ satisfying $\Red(\varphi)=(h,\alpha,\beta)$ if and only if for every positive integer $x$ all of the following are true:
\begin{enumerate}[label={\rm (O\arabic*)}]
\item\label{prop:odd:a} $\varphi^x(1)\equiv \pi_{\alpha}(x) \pmod{\ord(\alpha)}$;
\item\label{prop:odd:b} $\beta'(x) = \alpha(xm)/m$; and
\item\label{prop:odd:c} $\varphi^m(x)=\beta(x)$.
\end{enumerate}
\end{proposition}
\begin{proof}
First, we will show that if $\varphi$ is a skew morphism of $\ZZ_n$ and $\Red(\varphi)=(h,\alpha,\beta)$, then properties~\ref{prop:odd:a}, \ref{prop:odd:b} and \ref{prop:odd:c} are all true. Recall that $\alpha=\varphi'$, so by Lemma~\ref{lem:comp:e} it follows that the auto-order of $\varphi$ is $m$. Then, since the complexity of $\varphi$ has the opposite parity then the skew morphism $\alpha$ derived from it, we have $\beta=\varphi^m$. Next, by Proposition~\ref{prop:quo} we have $n/|\ker\varphi|=\ord(\varphi')=\ord(\alpha)$ and also $\varphi^x(1)\equiv \pi_{\alpha}(x) \pmod{n/|\ker\varphi|}$, and hence \ref{prop:odd:a} holds. Next, \ref{prop:odd:b} follows from the following equation that holds modulo $\ord(\beta)$. Note that we repeatedly use Lemma~\ref{lem:power}, Proposition~\ref{prop:quo}, and the definition of $\sigma\,$:
\begin{align*}
 m\beta'(x) &= m\sigma_{\beta}(x,1) = m\sigma_{\varphi^m}(x,1) = \sum_{0\leq i \leq x-1}m\pi_{\varphi^m}\left(\varphi^{mi}(1)\right) = \sum_{0\leq i \leq x-1} \sigma_{\varphi}(m,\varphi^{mi}(1)) \\
 &= \sigma_{\varphi}(mx,1) = \varphi'(mx) = \alpha(mx)\, .
\end{align*}
Finally, \ref{prop:odd:c} follows trivially from the fact that $\beta=\varphi^m$.

Next, we show that if properties~\ref{prop:odd:a}, \ref{prop:odd:b} and \ref{prop:odd:c} are all true, then the mapping $\varphi$ defined by Equation~\eqref{eq:odd} is a skew morphism of $\ZZ_n$ and $\Red(\varphi)=(h,\alpha,\beta)$. Note that \ref{prop:odd:c} implies that $\ord(\varphi)=\ord(\beta)m$. It follows that if $i$ is any element of $\ZZ_{\ord(\beta)m}$, then $\varphi^i$ is a well-defined permutation of $\ZZ_n$. As the first step, we use induction on $b$ to show that $\varphi(a+b)=\varphi(a)+\varphi^{\alpha^a(1)}(b)$ for any pair of non-zero elements $a,b$ in $\ZZ_n$. For $b=1$ we have $\varphi(a+1)=h + \beta^{k_1}(h) + \dots + \beta^{k_{a-1}}(h) + \beta^{k_{a}}(h) = \varphi(a) + \varphi^{mk_a}(\varphi(1))=\varphi(a) + \varphi^{\alpha^a(1)}(1)$. Next, since $\alpha^a(1)=(\alpha^a(1)-1)+1=mk_a+1$ and $\beta=\varphi^m$, we find that $\varphi^{\alpha^a(1)}(b)=\beta^{k_a}(\varphi(b))$. Then using the fact that $\varphi(b)=\varphi(b-1+1)=\varphi(b-1)+\varphi^{\alpha^{b-1}(1)}(1)$, Lemma~\ref{lem:powprod}, and the inductive hypothesis we find that 
\begin{align*}
\varphi(a)+\varphi^{\alpha^a(1)}(b) &= \varphi(a) + \beta^{k_a}\left(\varphi(b-1)+\varphi^{\alpha^{b-1}(1)}(1)\right) \\ 
&= \varphi(a) + \beta^{k_a}(\varphi(b-1)) +\beta^{\sigma_{\beta}(k_a,\varphi(b-1))} \left(\varphi^{\alpha^{b-1}(1)}(1)\right) \\
&= \varphi(a) + \varphi^{\alpha^a(1)}(b-1) + \beta^{\sigma_{\beta}(k_a,\varphi(b-1))}\left( \beta^{k_{b-1}}(\varphi(1)) \right) \\ 
&=  \varphi(a+b-1) + \beta^{\sigma_{\beta}(k_a,\varphi(b-1))}\left( \beta^{k_{b-1}}(h) \right)\, .
\end{align*}
On the other hand, we have $\varphi(a+b)=\varphi(a+b-1+1)=\varphi(a+b-1)+\varphi^{\alpha^{a+b-1}(1)}(1)=\varphi(a+b-1)+\beta^{k_{a+b-1}}(h)$. Hence to show that $\varphi(a+b)=\varphi(a)+\varphi^{\alpha^a(1)}(b)$ it is sufficient to verify that $\beta^{\sigma_{\beta}(k_a,\varphi(b-1))}\left( \beta^{k_{b-1}}(h) \right)$ is equal to $\beta^{k_{a+b-1}}(h)$ or, equivalently, that $\sigma_{\beta}(k_a,\varphi(b-1))+k_{b-1} \equiv k_{a+b-1} \pmod{\ord(\beta)}$. Next, note that 
\begin{align*}
k_{a+b-1} &= \frac{\alpha^{a+b-1}(1)-1}{m} = \frac{\alpha^{b-1}(\alpha^a(1))-1}{m}= \frac{\alpha^{b-1}(1+mk_a)-1}{m} 
= \frac{\alpha^{b-1}(1)+\alpha^{\sigma_{\alpha}(b-1,1)}(mk_a)-1}{m} \\ &= \frac{\alpha^{\sigma_{\alpha}(b-1,1)}(mk_a)+\alpha^{b-1}(1) - 1}{m} = \frac{\alpha^{\sigma_{\alpha}(b-1,1)}(mk_a)}{m} + k_{b-1}\, ,
\end{align*}
and so we need to prove that $\sigma_{\beta}(k_a,\varphi(b-1)) \equiv (\alpha^{\sigma_{\alpha}(b-1,1)}(mk_a))/m \pmod{\ord(\beta)}$. Next, since $\alpha$ has even complexity and auto-order $m$, by Theorem~\ref{thm:even} we deduce that $\alpha$ restricts to a skew morphism $\alpha_m$ of $\langle m\rangle$ given by $\alpha_m(x)=\alpha(xm)/m$. Also note that $\alpha$ is a skew morphism of $\ZZ_{\ord(\beta)m}$, and so $\alpha_m$ is a skew morphism of the cyclic group of order $\ord(\beta)$. Therefore \ref{prop:odd:b} can be rewritten as $\beta'(x)=\alpha_m(x)$, and so $\beta'$ and $\alpha_m$ are equal to the same skew morphism of $\ZZ_{\ord(\beta)}$. These observations and Proposition~\ref{prop:quo:a} imply that the congruence $\sigma_{\beta}(k_a,\varphi(b-1)) \equiv (\alpha^{\sigma_{\alpha}(b-1,1)}(mk_a))/m \pmod{\ord(\beta)}$ holds true if and only if $(\beta')^{\varphi(b-1)}(k_a) = (\beta')^{\sigma_{\alpha}(b-1,1)}(k_a)$, which is equivalent to $\varphi(b-1)\equiv \sigma_{\alpha}(b-1,1) \pmod{\ord(\beta')}$. Using \ref{prop:odd:c}, Equation~\ref{eq:odd}, and the fact that $h=\varphi(1)$, we find that $\varphi(b-1)=\varphi(1)+\varphi^{\alpha(1)}(1)+\dots+\varphi^{\alpha^{b-2}(1)}(1)$. Then, since $\beta'=\alpha_m$, we deduce that the order of $\beta'$ divides the order of $\alpha$, and hence from \ref{prop:odd:a} we obtain $\varphi(b-1)\equiv \pi_{\alpha}(1)+\pi_{\alpha}(\alpha(1))+\dots+\pi_{\alpha}(\alpha^{b-2}(1)) \pmod{\ord(\beta')}$. The right hand side of this congruence is by definition equal to $\sigma_{\alpha}(b-1,1)$, which proves our claim.  

We showed that $\varphi(a+b)=\varphi(a)+\varphi^{\alpha^a(1)}(b)$, and so if we let $\pi_{\varphi}(i)=\alpha^i(1)$, then $\varphi(a+b)=\varphi(a)+\varphi^{\pi_{\varphi}(a)}(b)$. (Recall that $\ord(\varphi)=\ord(\beta)m$ and $\alpha$ is a skew morphism of $\ZZ_{\ord(\beta)m}$, so this is a well-defined power function.) Next, we show that $\varphi(0)=0$. Recall that we defined $\varphi(0)$ as $\varphi(n)$, and this is equal to $\varphi(1)+\varphi^{\alpha(1)}(1)+\dots+\varphi^{\alpha^{n-1}(1)}(1)$. Then, since $\ord(\alpha)$ divides $n$, we deduce that $\varphi(n)=(\varphi(1)+\varphi^{\alpha(1)}(1)+\dots+\varphi^{\alpha^{\ord(\alpha)-1}(1)}(1))(n/\ord(\alpha))$. To show that the right hand side is divisible by $n$ (and consequently equal to $0$ in $\ZZ_n$) it is sufficient to prove that $\varphi(1)+\varphi^{\alpha(1)}(1)+\dots+\varphi^{\alpha^{\ord(\alpha)-1}(1)}(1)$ is a multiple of $\ord(\alpha)$. Note that by \ref{prop:odd:a} and Proposition~\ref{prop:quo:a} we have $\varphi(1)+\varphi^{\alpha(1)}(1)+\dots+\varphi^{\alpha^{\ord(\alpha)-1}(1)}(1) \equiv \pi_{\alpha}(1)+\pi_{\alpha}(\alpha(1))+\dots+\pi_{\alpha}(\alpha^{\ord(\alpha)-1}(1)) \equiv \sigma_{\alpha}(\ord(\alpha),1)\equiv \alpha'(\ord(\alpha)) \equiv \alpha'(0)\equiv 0 \pmod{\ord(\alpha)}$, and so we deduce that $\varphi(0)=0$. Finally, if $\varphi(a)=\varphi(b)$, then also $\varphi^{m-1}(\varphi(a))=\varphi^{m-1}(\varphi(b))$, and so $\beta(a)=\beta(b)$. But $\beta$ is a bijection, hence we must have $a=b$, and we deduce that $\varphi$ is a bijection too. It follows that $\varphi$ is a skew morphism with power function $\pi_{\varphi}$. Note that Equation~\eqref{eq:odd} implies that $\varphi(1)=h$, and from \ref{prop:odd:c} we deduce that $\beta=\varphi^m$. Hence to show that $\Red(\varphi)=(h,\alpha,\beta)$, it is sufficient to prove that $\varphi'=\alpha$. Since $\pi_{\varphi}(i)=\alpha^i(1)$, by Proposition~\ref{prop:quo:b} we deduce that $(\varphi')^i(1)=\alpha^i(1)$ for all positive integers $i$. In particular, by Theorem~\ref{thm:orderoforbit} it follows that $\ord(\varphi')=\ord(\alpha)$. Then by \ref{prop:odd:a} and Proposition~\ref{prop:quo:c} we find that $\pi_{\varphi'}(i)=\pi_{\alpha}(i)$, and so by Lemma~\ref{lem:skewfromorbit} we have $\varphi'=\alpha$.  
\end{proof}

If we want to find all skew morphisms of $\ZZ_n$ of order $\ell$ with odd complexity, we need to identify all candidates for the triple $(h,\alpha,\beta)$. By Proposition~\ref{prop:odd} and the prior discussion we know that $h$ is an element of $\ZZ_n$, and also that $\beta$ is a skew morphism of $\ZZ_n$ such that $\ord(\beta)$ divides $\ell$ and $\alpha$ is a skew morphism of $\ZZ_{\ell}$ with even complexity and auto-order $\ell/\ord(\beta)$. Also note that \ref{prop:odd:b} (in Proposition~\ref{prop:odd}) only refers to $\alpha$ and the skew morphism derived from $\beta$, so this relates $\alpha$ and $\beta$. Moreover, by Equation~\eqref{eq:odd} we have $\varphi(1)=h$, and hence from \ref{prop:odd:a} we deduce that $h\equiv \pi_{\alpha}(1) \pmod{\ord(\alpha)}$, which relates $h$ and $\alpha$. All this observations can be used to limit the number of all possible candidates for $(h,\alpha,\beta)$. Next, for each candidate we use Equation~\eqref{eq:odd} to define $\varphi$. Note that for each $a\geq 2$ we have $\varphi(a)=\varphi(a-1)+\beta^{k_{a-1}}(h)$, so the time complexity of this step is linear in $n$. As the final step we check if both \ref{prop:odd:a} and \ref{prop:odd:c} hold. In the first case we have $x\in \{1,2,\dots, \ord(\varphi)\}$ and in the second case we have $x\in \ZZ_n$. By Theorem~\ref{thm:orderofskew} we know that $\ord(\varphi)\leq n$, and so the time complexity of these checks is also linear in $n$. 

A quick computer search shows that if we omit any of the three properties \ref{prop:odd:a}, \ref{prop:odd:b} and \ref{prop:odd:c}, then there exists a triple $(h,\alpha,\beta)$ satisfying the two remaining properties such that Equation~\eqref{eq:odd} either does not define a skew morphism of $\ZZ_n$ or define a skew morphism $\varphi$ of $\ZZ_n$ but $\Red(\varphi)\neq (h,\alpha,\beta)$. It follows that none of the three properties is redundant.

\subsection{Skew morphisms with even complexity}%%%%%%  

Let $\varphi$ be a skew morphism of $\ZZ_n$ with even complexity and auto-order $m$, and let $\Red(\varphi)=(h,\alpha,\beta)$. Also let $\psi$ denote the restriction of $\varphi$ to the orbit of $\langle \varphi \rangle$ containing $1$. Then from Theorem~\ref{thm:even} we deduce that
\begin{align}\label{eq:even}
\psi^e(1)=1+fm+\beta^{\alpha(1)}(f)m+\dots + \beta^{\alpha(e-1)}(f)m \text{ for all } e\geq 2\, ,
\end{align}
and then by Lemma~\ref{lem:skewfromorbit} we have
\begin{equation}\label{eq:even2}
\varphi(a)=\psi(1)+\psi^{\alpha(1)}(1)+\psi^{\alpha^2(1)}(1)+\dots+\psi^{\alpha^{a-1}(1)}(1) \text{ for all } a\in \ZZ_n\, . 
\end{equation}
In the following proposition we characterise all triples $(h,\alpha,\beta)$ for which Equations~\eqref{eq:even} and~\eqref{eq:even2} define a skew morphism of $\ZZ_n$. Similar to the previous case we define $\varphi(0)$ as $\varphi(n)$. 

\begin{proposition}\label{prop:even}
Let $n$ be a positive integer, let $m$ be a divisor of $n$, let $\ell$ be a positive integer strictly smaller than $n$, let $h$ be an element of $1+\langle m\rangle$, let $f=(h-1)/m$, let $\beta$ be a skew morphism of $\ZZ_{n/m}$ such that $\ord(\beta)$ divides $\ell$, and let $\alpha$ be a skew morphism of $\ZZ_{\ell}$ with odd complexity and auto-order $m$ such that $\ord(\alpha)$ divides $n$. Also let $\psi$ and $\varphi$ be the mappings given by Equations~\eqref{eq:even} and~\eqref{eq:even2}. Then $\varphi$ is a skew morphism of $\ZZ_n$ satisfying $\Red(\varphi)=(h,\alpha,\beta)$ if and only if for every positive integer $x$ all of the following are true:
\begin{enumerate}[label={\rm (E\arabic*)}]
\item\label{prop:even:a} $\psi^x(1)\equiv \pi_{\alpha}(x) \pmod{\ord(\alpha)}$;
\item\label{prop:even:b} $\beta'(x) \equiv \alpha^m(x) \pmod{\ord(\beta)}$;
\item\label{prop:even:c} $\varphi(xm)=\beta(x)m$;
\item\label{prop:even:d} $\varphi(x+m)=\varphi(x)+\varphi^{\alpha^x(1)}(m)$; and
\item\label{prop:even:e} the smallest positive integer $i$ such that $\psi^i(1)=1$ is equal to $\ell$.
\end{enumerate}
\end{proposition}
\begin{proof}
First, we will show that if $\varphi$ is a skew morphism of $\ZZ_n$ and $\Red(\varphi)=(h,\alpha,\beta)$, then properties \ref{prop:even:a}, \ref{prop:even:b}, \ref{prop:even:c}, \ref{prop:even:d} and \ref{prop:even:e} are all true. Since $\alpha=\varphi'$, by Lemma~\ref{lem:comp:e} it follows that the auto-order of $\varphi$ is $m$. Then, since $\alpha$ has odd complexity, we deduce that $\varphi$ has even complexity, hence $\beta=\varphi\vert_{\langle m\rangle}$. Then it can be easily checked (see the proof of Theorem~\ref{thm:even}) that $\varphi^x(1)=\psi^x(1)$, and then the same argument as in the proof of Proposition~\ref{prop:odd} implies \ref{prop:even:a}. Next, \ref{prop:even:b} follows from the following equation that holds modulo $\ord(\beta)$. Note that we use Proposition~\ref{prop:quo:a}, the definition of $\sigma$, and the fact that for an arbitrary exponent $e$ we have $\pi_{\beta}(\beta^e(1))\equiv \pi_{\varphi}(\varphi^e(m)) \pmod{\ord(\beta)}$ (this is a consequence of the fact that $\beta$ is the restriction of $\varphi$ to $\langle m\rangle$):
\begin{align*}
\alpha^m(x) = (\varphi')^m(x) = \sigma_{\varphi}(x,m) = \sum_{0\leq i \leq x-1}\pi_{\varphi}(\varphi^i(m)) = \sum_{0\leq i \leq x-1}\pi_{\beta}(\beta^i(1))=\sigma_{\beta}(x,1)=\beta'(x)\, .
\end{align*}
Further, \ref{prop:even:c} follows from the fact that $\beta=\varphi\vert_{\langle m\rangle}$, \ref{prop:even:d} is a straightforward consequence of the definition of a skew morphism and Proposition~\ref{prop:quo:b}. Finally, since $\psi$ defines the orbit of $\langle \varphi\rangle$ that contains $1$, by Theorem~\ref{thm:orderoforbit} we know that the smallest positive integer $i$ such that $\psi^i(1)=1$ is the order of $\varphi$. Recall that $\varphi'$ is a skew morphism of $\ZZ_{\ord(\varphi)}$, and so $i=\ord(\varphi)=\ell$. This proves \ref{prop:even:e}.

We proceed to show that if properties \ref{prop:even:a}, \ref{prop:even:b}, \ref{prop:even:c}, \ref{prop:even:d} and \ref{prop:even:e} are all true, then the mapping $\varphi$ defined by Equations~\eqref{eq:even} and~\eqref{eq:even2} is a skew morphism of $\ZZ_n$ and $\Red(\varphi)=(h,\alpha,\beta)$. Throughout the rest of the proof we will often use elements of $\ZZ_{\ell}$ (and their images under $\alpha$) as exponents of various permutations. Note that for this to be well-defined, the length of the cycle (of the given permutation on the given element) must divide $\ell$. For example, by \ref{prop:even:e} we know that for any $i\in \ZZ_{\ell}$ we can write $\psi^i(1)$, and since $\ord(\beta)$ divides $\ell$ we can write $\beta^i(x)$ for every $x\in \ZZ_{n/m}$. Moreover, by \ref{prop:even:c} we deduce that we also have $\varphi^{\ell}(xm)=\beta^{\ell}(x)m=xm$, and so we can also write $\varphi^i(xm)$. As the first step, we will show that for each positive integer $y$ we have $\varphi(x+ym)=\varphi(x)+\varphi^{\alpha^x(1)}(ym)$. To prove this we use induction on $y$. For $y=1$ this is equivalent to \ref{prop:even:d}, so the claim is trivially true. Next, we apply the inductive hypothesis to $y-1$ and deduce that $\varphi(x+ym)=\varphi(x+m+(y-1)m)=\varphi(x+m)+\varphi^{\alpha^{x+m}(1)}((y-1)m)$. Using \ref{prop:even:d} and \ref{prop:even:c} this can be further rewritten as $\varphi(x)+\varphi^{\alpha^x(1)}(m)+\varphi^{\alpha^{x+m}(1)}((y-1)m)=\varphi(x)+\left(\beta^{\alpha^x(1)}(1)+ \beta^{\alpha^{x+m}(1)}(y-1)\right)m$. The term in parentheses can be further simplified using \ref{prop:even:b}, Lemma~\ref{lem:derived}, and Lemma~\ref{lem:powprod}. Note that we use the fact that $\beta$ is a skew morphism:
\begin{align*}
\beta^{\alpha^x(1)}(1)+ \beta^{\alpha^{x+m}(1)}(y-1) &= \beta^{\alpha^x(1)}(1)+ \beta^{\alpha^{m}(\alpha^x(1))}(y-1)=\beta^{\alpha^x(1)}(1)+ \beta^{\beta'(\alpha^x(1))}(y-1) \\
&= \beta^{\alpha^x(1)}(1)+ \beta^{\sigma_{\beta}(\alpha^x(1),1)}(y-1)=\beta^{\alpha^x(1)}(1+y-1)=\beta^{\alpha^x(1)}(y)\, .
\end{align*}
Hence we have $\varphi(x+ym)= \varphi(x)+\beta^{\alpha^x(1)}(y)m=\varphi(x)+\varphi^{\alpha^x(1)}(ym)$.

Next, we use induction to prove that for each positive integer $e$ we have $\varphi^e(1)=\psi^e(1)$. For $e=1$ this follows from Equation~\eqref{eq:even2} (with $a=1$). Note that Equation~\eqref{eq:even} implies that $\psi^i(1)=\psi^{i-1}(1)+\beta^{\alpha(i-1)}(f)m$ for each positive integer $i$. Hence, using the inductive hypothesis we find that $\varphi^e(1)=\varphi(\varphi^{e-1}(1))=\varphi(\psi^{e-1}(1))=\varphi(\psi^{e-2}(1)+\beta^{\alpha(e-2)}(f)m)$. Next, we use the identity proved in the previous paragraph to find that $\varphi^e(1) = \varphi(\psi^{e-2}(1)) + \varphi^{\alpha^{\psi^{e-2}(1)}(1)}(\beta^{\alpha(e-2)}(f)m)$. Using the inductive hypothesis twice (on the first term), and \ref{prop:even:a} and \ref{prop:even:c} (on the second term) this can be rewritten as $\psi^{e-1}(1)+\varphi^{\alpha^{\pi_{\alpha}(e-2)}(1)}(\varphi^{\alpha(e-2)}(fm))=\psi^{e-1}(1)+\varphi^{\alpha^{\pi_{\alpha}(e-2)}(1)+\alpha(e-2)}(fm)$. Since $\alpha$ is a skew morphism we have $\alpha^{\pi_{\alpha}(e-2)}(1)+\alpha(e-2)=\alpha(e-2)+\alpha^{\pi_{\alpha}(e-2)}(1)=\alpha(e-2+1)=\alpha(e-1)$, and so we deduce that $\varphi^e(1)=\psi^{e-1}(1)+\varphi^{\alpha(e-1)}(fm)=\psi^{e-1}(1)+\beta^{\alpha(e-1)}(f)m=\psi^e(1)$, which proves our claim. Note that from now on we can write $\varphi^i(1)$ for any $i\in \ZZ_{\ell}$.  

As the next step we show that $\varphi(x+1+ym)=\varphi(x)+\varphi^{\alpha^x(1)}(1+ym)$. We first need to prove that $\varphi^{\alpha^x(1)}(1+ym)$ is well-defined, and so we need to prove that $\varphi^{\ell}(1+ym)=1+ym$. Let $e$ be an arbitrary positive integer and note that $\varphi^{e}(1+ym) = \varphi^{e-1}(\varphi(1+ym))=\varphi^{e-1}(\varphi(1)+\varphi^{\alpha(1)}(ym))$. Then, since \ref{prop:even:c} implies that $\varphi^{\alpha(1)}(ym)$ is a multiple of $m$, we find that $\varphi^{e-1}(\varphi(1)+\varphi^{\alpha(1)}(ym))=\varphi^{e-2}(\varphi^2(1)+\varphi^{\alpha(1)+\alpha^{\varphi(1)}(1)}(ym))$. The second term will always be a multiple of $m$, and so by repeated application of $\varphi$ we find that $\varphi^{e}(1+ym)=\varphi^{e}(1)+\varphi^{\alpha(1)+\alpha^{\varphi(1)}(1)+\dots+\alpha^{\varphi^{e-1}(1)}(1)}(ym)$. Since $\varphi^e(1)\equiv \psi^e(1)\equiv \pi_{\alpha}(e) \pmod{\ord(\alpha)}$ for each exponent $e$, this can be further rewritten as $\varphi^{e}(1)+\varphi^{\alpha(1)+\alpha^{\pi_{\alpha}(1)}(1)+\dots+\alpha^{\pi_{\alpha}(e-1)}(1)}(ym)$. Lemma~\ref{lem:skewfromorbit} can be used to simplify the second exponent and we obtain $\varphi^{e}(1+ym)=\varphi^{e}(1)+\varphi^{\alpha(e)}(ym)$. It follows that $\varphi^{\ell}(1+ym)=\varphi^{\ell}(1)+\varphi^{\alpha(\ell)}(ym)=1+\varphi^{\ell}(ym)=1+ym$, and so $\varphi^{\alpha^x(1)}(1+ym)$ is well-defined. The equation $\varphi(x+1+ym)=\varphi(x)+\varphi^{\alpha^x(1)}(1+ym)$ now follows easily as $\varphi(x)+\varphi^{\alpha^x(1)}(1+ym)=\varphi(x)+\varphi^{\alpha^x(1)}(1)+\varphi^{\alpha(\alpha^x(1))}(ym)=\varphi(x)+\varphi^{\alpha^x(1)}(1)+\varphi^{\alpha^{x+1}(1)}(ym)=\varphi(x+1)+\varphi^{\alpha^{x+1}(1)}(ym)=\varphi(x+1+ym)$.

We are now ready to prove that for all $x,y\in\ZZ_n$ we have $\varphi(x+y)=\varphi(x)+\varphi^{\alpha^x(1)}(y)$ and $\varphi^{\ell}(y)=y$. We use induction on $y$. For $y=1$ we have $\varphi^{\ell}(1)=\psi^{\ell}(1)=1$ and $\varphi(x+1)=\psi(1)+\psi^{\alpha(1)}(1)+\dots+\psi^{\alpha^{x-1}(1)}(1)+\psi^{\alpha^{x}(1)}(1)=\varphi(x)+\psi^{\alpha^{x}(1)}(1)=\varphi(x)+\varphi^{\alpha^{x}(1)}(1)$, and so the claim is true. Next, we will show that $\varphi^{\ell}(y+1)=y+1$ and $\varphi(x+y+1)=\varphi(x)+\varphi^{\alpha^x(1)}(y+1)$. Let $e$ be an arbitrary positive integer. We first show that $\varphi^{e}(y+1)=\varphi^e(y)+\varphi^{\alpha^y(e)}(1)$. The left hand side of the equation can be rewritten as $\varphi^{e-1}(\varphi(y+1))=\varphi^{e-1}(\varphi(y)+\varphi^{\alpha^y(1)}(1))$. Recall that $\varphi^{\alpha^y(1)}(1)$ is equal to $\psi^{\alpha^y(1)}(1)$, which by Equation~\eqref{eq:even} has the form $1+zm$ for some integer $z$. Hence we can use the identity from the previous paragraph to deduce that $\varphi^{e-1}(\varphi(y)+\varphi^{\alpha^y(1)}(1))=\varphi^{e-2}(\varphi(\varphi(y)+\varphi^{\alpha^y(1)}(1)))=\varphi^{e-2}(\varphi^2(y)+\varphi^{\alpha^y(1)+\alpha^{\varphi(y)}(1)}(1))$. The second term is again of the form $1+zm$, and by repeated application of $\varphi$ we find that $\varphi^{e}(y+1)=\varphi^{e}(y)+\varphi^{\alpha^y(1)+\alpha^{\varphi(y)}(1)+\dots+\alpha^{\varphi^{e-1}(y)}}(1)$. We want to prove that this is equal to $\varphi^e(y)+\varphi^{\alpha^y(e)}(1)$, and so we need to show that $\alpha^{y}(e)\equiv \alpha^y(1)+\alpha^{\varphi(y)}(1)+\dots+\alpha^{\varphi^{e-1}(y)} \pmod{\ell}$. To prove this, first note that by Equation~\eqref{eq:even2} and \ref{prop:even:a} we have $\varphi(z)\equiv \pi_{\alpha}(1)+\pi_{\alpha}(\alpha(1))+\dots+\pi_{\alpha}(\alpha^{z-1}(1))\equiv\sigma_{\alpha}(z,1) \pmod{\ord(\alpha)}$ for each $z\in \ZZ_n$. It follows that 
\begin{align*}
\alpha^{y}(e)&=\alpha^{y}(1+e-1) = \alpha^y(1)+\alpha^{\sigma_{\alpha}(y,1)}(e-1) = \alpha^y(1)+\alpha^{\varphi(y)}(e-1) \\
&= \alpha^y(1)+\alpha^{\varphi(y)}(1+e-2) = \alpha^y(1)+\alpha^{\varphi(y)}(1)+\alpha^{\sigma_{\alpha}(\varphi(y),1)}(e-2) \\
&= \alpha^y(1)+\alpha^{\varphi(y)}(1)+\alpha^{\varphi^2(y)}(e-2) = \cdots = \alpha^y(1)+\alpha^{\varphi(y)}(1)+\dots+\alpha^{\varphi^{e-1}(y)}(1)\, ,
\end{align*}  
which proves our claim. In particular, if we let $e=\ell$, then $\varphi^{\ell}(y+1)=\varphi^{\ell}(y)+\varphi^{\alpha^y(\ell)}(1)=\varphi^{\ell}(y)+\varphi^{\ell}(1)=\varphi^{\ell}(y)+1$. We may now apply the inductive hypothesis to deduce that $\varphi^{\ell}(y)=y$, and so $\varphi^{\ell}(y+1)=y+1$. We proceed to show that $\varphi(x+y+1)=\varphi(x)+\varphi^{\alpha^x(1)}(y+1)$. Note that by the inductive hypothesis we have $\varphi(x+y+1)=\varphi(x+y)+\varphi^{\alpha^{x+y}(1)}(1)=\varphi(x)+\varphi^{\alpha^x(1)}(y)+\varphi^{\alpha^{x+y}(1)}(1)$, so we need to show that $\varphi^{\alpha^x(1)}(y)+\varphi^{\alpha^{x+y}(1)}(1)=\varphi^{\alpha^x(1)}(y+1)$. But this is an immediate consequence of the identity $\varphi^e(y)+\varphi^{\alpha^y(e)}(1)=\varphi^{e}(y+1)$ with $e=\alpha^x(1)$. This proves that $\varphi(x+y)=\varphi(x)+\varphi^{\alpha^x(1)}(y)$.

Similar to the proof of Proposition~\ref{prop:odd}, we now define $\pi_{\varphi}(i)=\alpha^i(1)$, and by the same argument we deduce that $\varphi(a+b)=\varphi(a)+\varphi^{\pi_{\varphi(a)}}(b)$ and $\varphi(0)=0$. To show that $\varphi$ is bijection suppose to the contrary that there exists a pair of elements $x$ and $y$ of $\ZZ_n$ such that $\varphi(x)=\varphi(x+y)$ and $y\neq 0$. Then $\varphi(x)=\varphi(x)+\varphi^{\pi_{\varphi}(x)}(y)$, and so $\varphi^{\pi_{\varphi}(x)}(y)=0$ for a non-zero element $y$ of $\ZZ_n$. It follows that there exists some non-zero element $z$ of $\ZZ_n$ such that $\varphi(z)=0$. Then from \ref{prop:even:d} we deduce that $\varphi(z+m)=\varphi^{\pi_{\alpha}(z)}(m)$. Note that Equations~\eqref{eq:even} and~\eqref{eq:even2} imply that $\varphi(a)$ is a multiple of $m$ if and only if $a$ is a multiple of $m$, and so we deduce that $z+m$, and consequently also $z$ is a multiple of $m$. In particular, we have $z=z'm$ for some non-zero element $z'$ of $\ZZ_{n/m}$. But then \ref{prop:even:c} implies that $0=\varphi(z)=\varphi(z'm)=\beta(z')m$, and so $\beta(z')$ must be the zero element of $\ZZ_{n/m}$, contradiction. It follows that $\varphi$ is a permutation, and consequently a skew morphism of $\ZZ_n$.

As the last step we show that $\Red(\varphi)=(h,\alpha,\beta)$. Equation~\eqref{eq:even} implies that $\varphi(1)=1+fm=h$, and from \ref{prop:even:c} we deduce that $\beta=\varphi\vert_{\langle m\rangle}$. Finally, since $\varphi^x(1)\equiv\psi^x(1)\equiv\pi_{\alpha}(x) \pmod{\ord(\alpha)}$, the same argument as in the proof of Proposition~\ref{prop:odd} yields $\varphi'=\alpha$.  
\end{proof}

If we want to find all skew morphisms of $\ZZ_n$ of order $\ell$ with even complexity, we need to identify all candidates for the triple $(h,\alpha,\beta)$. By Proposition~\ref{prop:even} and the prior discussion we know that $h$ is an element of $\ZZ_n$, and also that $\alpha$ is a skew morphism of $\ZZ_{\ell}$ with odd complexity and $\beta$ is a skew morphism of $\ZZ_{n/m}$ where $m$ is the auto-order of $\alpha$. Also note that \ref{prop:even:b} (in Proposition~\ref{prop:even}) only refers to $\alpha$ and the skew morphism derived from $\beta$, so this relates $\alpha$ and $\beta$. Moreover, by Equations~\eqref{eq:even} and ~\eqref{eq:even2} we have $\varphi(1)=h$, and hence from \ref{prop:even:a} we deduce that $h\equiv \pi_{\alpha}(1) \pmod{\ord(\alpha)}$, which relates $h$ and $\alpha$. All this observations can be used to limit the number of all possible candidates for $(h,\alpha,\beta)$. Next, for each candidate we use Equations~\eqref{eq:even} and ~\eqref{eq:even2} to define $\varphi$. Note that for each $e\geq 2$ we have $\psi^e(1)=\psi^{e-1}(1)+\beta^{\alpha(e-1)}(f)m$ and for each $a\geq 2$ we have $\varphi(a)=\varphi(a-1)+\psi^{\alpha^{a-1}(1)}(1)$, so the time complexity of this step is linear in $n$. As the final step we check if properties \ref{prop:even:a}, \ref{prop:even:c}, \ref{prop:even:d} and \ref{prop:even:e} all hold. In the first case we have $x\in \{1,2,\dots, \ord(\varphi)\}$, in the second and third case we have $x\in \ZZ_n$, and in the fourth case we have $\ell=n/|\ker\varphi|\leq n$. By Theorem~\ref{thm:orderofskew} we know that $\ord(\varphi)\leq n$, and so the time complexity of these checks is also linear in $n$. 

A computer search shows that if we omit any of the properties \ref{prop:even:a}, \ref{prop:even:c}, \ref{prop:even:d} and \ref{prop:even:e}, then there exists a triple $(h,\alpha,\beta)$ satisfying \ref{prop:even:b} and the three remaining properties such that Equations~\eqref{eq:even} and ~\eqref{eq:even2} either does not define a skew morphism of $\ZZ_n$ or define a skew morphism $\varphi$ of $\ZZ_n$ but $\Red(\varphi)\neq(h,\alpha,\beta)$. In contrast, it can be shown that \ref{prop:even:b} follows from the other four properties listed in Proposition~\ref{prop:even}, making it redundant. Nevertheless, we include this property in the statement of Proposition~\ref{prop:even} as it provides a useful constrain linking $\alpha$ and $\beta$, which can be used to considerably reduce the search space of triples $(h,\alpha,\beta)$ that may define a skew morphism of $\ZZ_n$. This makes \ref{prop:even:b} useful for both theoretical and practical applications of Proposition~\ref{prop:even}. Here we provide a proof that \ref{prop:even:b} follows from the other four properties:
\begin{proof}
Let $h$, $\alpha$, $\beta$, $\psi$ and $\varphi$ be defined as in Proposition~\ref{prop:even}, and suppose that \ref{prop:even:a}, \ref{prop:even:c}, \ref{prop:even:d} and \ref{prop:even:e} are all true. We will show that $\beta'(x) \equiv \alpha^m(x) \pmod{\ord(\beta)}$. Note that $\beta'$ is a skew morphism of $\ZZ_{\ord(\beta)}$, and so the values $\beta'(x)$ taken modulo $\ord(\beta)$ are exactly same as the values $\beta'(x)$. As $\alpha$ has odd complexity and auto-order $m$, we deduce that $\alpha^m$ is a well-defined skew morphism of $\ZZ_{\ell}$. Then since $\ord(\beta)$ divides $\ell$, we deduce that the values $\alpha^m(x)$ taken modulo $\ord(\beta)$ are well-defined. As the first step we will show that for each positive integer $e$ we have $(\beta')^e(1)\equiv (\alpha^m)^e(1) \pmod{\ord(\beta)}$. To prove this first note that by \ref{prop:even:d} and \ref{prop:even:c} we have the following: 
\[  \varphi(em+m) = \varphi(em)+\varphi^{\alpha^{em}(1)}(m)= \left(\beta(e)+\beta^{\alpha^{em}(1)}(1) \right)m=\left(\beta(e)+\beta^{(\alpha^m)^e(1)}(1) \right)m \, . \]
On the other hand, using $\ref{prop:even:c}$ and Proposition~\ref{prop:quo:b} we deduce that 
\[  \varphi(em+m) = \beta(e+1)m= \left(\beta(e)+\beta^{\pi_{\beta}(e)}(1) \right)m=\left(\beta(e)+\beta^{(\beta')^e(1)}(1) \right)m \, . \]
Since $\beta$ is a skew morphism of $\ZZ_{n/m}$, it follows that $\beta(e)+\beta^{(\alpha^m)^e(1)}(1)$ is equal to $\beta(e)+\beta^{(\beta')^e(1)}(1)$, therefore $(\beta')^e(1) \equiv (\alpha^m)^e(1) \pmod{\ord(\beta)}$ as claimed. Note that by Theorem~\ref{thm:orderoforbit} this implies that $\ord(\beta')$ divides $\ord(\alpha^m)$, and consequently also $\ord(\alpha)$.

Next, we show that for each $x\in \ZZ_{\ord(\beta')}$ we have $\pi_{\beta'}(x)\equiv \pi_{\alpha^m}(x) \pmod{\ord(\beta')}$. Using Proposition~\ref{prop:quo:b}, the definition of $\varphi^*$, and $\ref{prop:even:c}$ we find that $\pi_{\beta'}(x)\equiv (\beta^*)^x(1)\equiv \beta^x(1) \equiv \varphi^x(m)/m \pmod{\ord(\beta')}$. Now we use induction on $x$ to show that also $\pi_{\alpha^m}(x)\equiv \varphi^x(m)/m \pmod{\ord(\beta')}$. (We will actually prove that this congruence holds modulo $\ord(\alpha)$, which is sufficient as $\ord(\beta')$ divides $\ord(\alpha)$.) For $x=1$ we use Lemma~\ref{lem:power} to deduce that $\pi_{\alpha^m}(1)\equiv \sigma_{\alpha}(m,1)/m \pmod{\ord(\alpha)}$. Using the definition of $\sigma_{\alpha}$ and \ref{prop:even:a} we find that $\sigma_{\alpha}(m,1)\equiv \psi(1)+\psi^{\alpha(1)}(1)+\dots+\psi^{\alpha^{m-1}(1)}(1)  \pmod{\ord(\alpha)}$, and the right hand side of this congruence is by definition equal to $\varphi(m)$. It follows that $\pi_{\alpha^m}(1)\equiv \varphi(m)/m \pmod{\ord(\alpha)}$. Next, we will show that if the claim is true for $x$, then also $\pi_{\alpha^m}(x+1)\equiv \varphi^{x+1}(m)/m \pmod{\ord(\beta')}$. From Lemma~\ref{lem:power} we deduce that $\pi_{\alpha^m}(x+1)\equiv \sigma_{\alpha}(m,x+1)/m \pmod{\ord(\alpha)}$. By Proposition~\ref{prop:quo:a} we find that $\sigma_{\alpha}(m,x+1)=(\alpha')^{x+1}(m)=\alpha'((\alpha')^x(m))=\alpha'(\sigma_{\alpha}(m,x))$. Using the inductive hypothesis and Lemma~\ref{lem:power} this can be further rewritten as $\alpha'(\pi_{\alpha^m}(x)m)=\alpha'(\varphi^x(m))=\sigma_{\alpha}(\varphi^x(m),1)$. Finally, by the definition of $\sigma_{\alpha}$ and \ref{prop:even:a} this is congruent to $\varphi(\varphi^x(m))$ modulo $\ord(\alpha)$, which proves our claim.

We have shown that $\alpha^m$ taken modulo $\ord(\beta)$ and $\beta'$ has the same orbit that contains $1$, and also that $\alpha^m$ and $\beta'$ have the same power function when taken modulo $\ord(\beta')$. It is now straightforward to check that $\alpha^m$ taken modulo $\ord(\beta)$ is a well-defined skew morphism of $\ZZ_{\ord(\beta)}$, and by Lemma~\ref{lem:skewfromorbit} we deduce that this skew morphism is equal to $\beta'$. It follows that $\ref{prop:even:b}$ is true.
\end{proof}  

\subsection{Database of skew morphisms of cyclic groups}%%%%%%  

Based on the observations from the previous sections, we have developed an algorithm for finding skew morphisms of finite cyclic groups. A \CC\ implementation of this algorithm is available at~\cite{alg}, and the corresponding online database of skew morphisms of finite cyclic groups can be found at~\cite{census}. While the algorithm is built on the ideas and techniques presented in this paper, it is regularly updated. As a result, some components may incorporate newer, unpublished observations. 

%%%%%%%%%%%
\section{Applications of the recursive characterisation}\label{sec:apps}
%%%%%%%%%%%

In this section we present various applications of the recursive characterisation of skew morphisms of cyclic groups. While all of the following findings are theoretical, we repeatedly used~\cite{census} for both proposing and testing various hypotheses.
 
\subsection{Complexities of skew morphisms}%%%%%% 

In what follows let $\Comp(\ZZ_n)$ denote the set of complexities of all skew morphisms of $\ZZ_n$. Note that $\ZZ_n$ admits a trivial automorphism for each positive integer $n$, and a non-trivial automorphism for each $n\geq 3$. It follows that $\{0,1\}\subseteq \Comp(\ZZ_n)$ for every $n\geq 3$. The equality $\{0,1\} = \Comp(\ZZ_n)$ holds if and only if $n\geq 3$ and $\ZZ_n$ admits no proper skew morphism, which is known to be true only if $n=4$ or $\gcd(n,\phi(n))=1$; see~\cite{KovacsNedela2011}. Skew morphisms with complexity at most $2$ are known as coset-preserving skew morphisms or smooth skew morphisms. All such skew morphisms for cyclic groups were classified in~\cite{BachratyJajcay2017}, and recently in~\cite{Bachraty} (and also independently in~\cite{HuKovacsKwon2024}) it was shown that all skew morphisms of $\ZZ_n$ have complexity at most $2$ if and only if $n=2^em$ with $e\in\{0,1,2,3,4\}$ and $m$ odd and square-free. A next natural step is to determine all $n$ such that all skew morphisms of $\ZZ_n$ have complexity at most $3$. We leave this as an open question, but throughout this section we provide some useful observations. 

\begin{proposition}\label{prop:psquare}
If $p$ is an odd prime, then $\Comp(\ZZ_{p^2})=\{0,1,3\}$.
\end{proposition}
\begin{proof}
First note that clearly $\{0,1\} \subseteq \Comp(\ZZ_{p^2})$. Next, let $\varphi$ be a proper skew morphism of $\ZZ_{p^2}$. Recall that the kernel of a skew morphism of a non-trivial group is always non-trivial, and so $|\ker\varphi|$ is either $p$ or $p^2$, but the latter case can be excluded as $\varphi$ is proper. It follows that $|\ker\varphi|=p$, and from Lemma~\ref{lem:star} we deduce that $\varphi^*$ is a skew morphism of $\ZZ_p$. But $\ZZ_p$ has no proper skew morphisms, and so the complexity of $\varphi^*$ is at most $1$. Furthermore, by Proposition~\ref{prop:quo:e} we have $\varphi^*=\varphi^{(2)}$, and so we deduce that the complexity of $\varphi$ is at most $3$. In particular, the complexity of a proper skew morphism of $\ZZ_{p^2}$ can be only $2$ or $3$. In what follows we will exclude the former case. 

Suppose to the contrary that there exists a skew morphism of $\ZZ_{p^2}$ with complexity $2$. By Lemma~\ref{lem:comp:d} we know that $\varphi$ is proper, and so our prior arguments imply that $\varphi^*$ is a skew morphism of $\ZZ_p$. Moreover, the complexity of $\varphi^*$ is $0$, hence Lemma~\ref{lem:comp:b} implies that it is the trivial permutation of $\ZZ_p$, and so the auto-order of $\varphi$ is $p$. Then by Lemma~\ref{lem:main} we find that $\varphi$ taken modulo $p$ is the trivial permutation of $\ZZ_{p^2}$, and it follows easily that $\ord(\varphi)\leq p$. On the other hand, by Theorem~\ref{thm:orderofskewcyclic} we have $\gcd(\ord(\varphi),p^2)>1$, hence $\ord(\varphi)=p$. It follows that $\varphi'$ is a skew morphism of $\ZZ_p$, and from Lemma~\ref{lem:comp:f} we know that the order of $\varphi'$ should be $p$, but by Theorem~\ref{thm:orderofskew} this cannot happen. Hence there are no skew morphisms of $\ZZ_{p^2}$ with complexity $2$, and so every proper skew morphism of $\ZZ_{p^2}$ has complexity $3$. The proof of existence of a proper skew morphism of $\ZZ_{p^2}$  (for all odd primes) can be found in~\cite[Proposition~4.9]{KovacsNedela2011}.  
\end{proof}

Another natural question about complexities of skew morphisms is whether there exist skew morphisms (of finite cyclic groups) of arbitrarily large complexities. The following theorem shows that the answer is positive.

\begin{theorem}\label{thm:unbounded}
Let $p$ be an odd prime, and let $e$ be a non-negative integer. Then the largest complexity of a skew morphism of $\ZZ_{p^e}$ is $2e-1$.  
\end{theorem}   
\begin{proof}
First, we use induction on $e$ to show that if $\varphi$ is a skew morphism of $\ZZ_{p^e}$ with complexity $c$, then $c\leq 2e-1$. Recall that every skew morphism of $\ZZ_p$ must be an automorphism of $\ZZ_p$, and so its complexity is at most $1$. Next, let $\varphi$ be a non-trivial skew morphism of $\ZZ_{p^e}$ with complexity $c$. Since $\ker\varphi$ is a non-trivial subgroup of $\ZZ_{p^e}$, it follows that $\varphi^*$ is a skew morphism of $\ZZ_{p^f}$ for some $f\leq e-1$. We may now apply the inductive hypothesis to deduce that the complexity of $\varphi^*$ is at most $2(e-1)+1$, and the rest follows from the fact that the complexity of $\varphi^*$ is $c-2$.

Next, we will show that there exists a skew morphism of $\ZZ_{p^e}$ with complexity $2e-1$. We will use the classification of skew morphisms of cyclic $p$-groups presented in~\cite{KovacsNedela2017}. Namely, if $e\geq 2$, then skew morphisms of $\ZZ_{p^e}$ are exactly the permutations $s_{i,j,k,l}=b_j^{\, -1}a^ib^kb_lb_j$, where $a$ is the automorphism of $\ZZ_{p^e}$ given by $x\mapsto (p+1)x$, $b$ is an automorphism of $\ZZ_{p^e}$ of order $p-1$, $b_j$ is the permutation of $\ZZ_{p^e}$ defined as $b_j(x)=1+(p+1)^j+\dots+(p+1)^{(x-1)j}$, and the integers $i$, $j$, $k$, $l$ satisfy the following conditions:
\begin{enumerate}[label={\rm (C\arabic*)}]
\item\label{cond:1} $i,l\in\{0,\dots,p^{e-1}-1 \}$, $k\in \{0,\dots,p-2\}$, $j\in \{0,\dots,p^{e-2-c}-1\}$, where $p^c=\gcd(i,p^{e-2})$;
\item\label{cond:2} if $i=0$ or $k=0$, then $l=0$; and
\item\label{cond:3} if $i\neq 0$ and $k\neq 0$, then $p^c\mid j$ and $p^{\max\{c,e-2-c\}}\mid l$.
\end{enumerate}
It can be easily verified that the $4$-tuple $(1,1,1,p^{e-2})$ satisfies all three conditions \ref{cond:1}, \ref{cond:2} and \ref{cond:3}, and we will show that the skew morphism $s_{1,1,1,p^{e-2}}$ defined as $b_1^{\, -1}abb_{p^{e-2}}b_1$ has complexity $2e-1$. Note that for any $j$ we have $b_j(1)=1$. Moreover, it can be easily seen that $(p+1)^{p^{e-2}}\equiv p^{e-1}+1 \pmod{p^e}$, and so $b_{p^{e-2}}$ is given by $b_{p^{e-2}}(x)=1+(p^{e-1}+1)+\dots+(p^{e-1}+1)^{x-1}$. Using the binomial theorem (and the fact that $b_{p^{e-2}}$ is a permutation of $\ZZ_{p^e}$ where $p$ is an odd prime) this can be further rewritten as $b_{p^{e-2}}(x)=1+(p^{e-1}+1)+(2p^{e-1}+1)\dots+((x-1)p^{e-1}+1)=x+p^{e-1}x(x-1)2^{-1}$. In particular, if $x$ is congruent to $0$ or $1$ modulo $p$, then $b_{p^{e-2}}(x)=x$. Next, we show that for each $x\in \ZZ_{p^e}$ the largest power of $p$ that divides $x$ is equal to the largest power of $p$ that divides $b_1(x)$. Let $p^f$ be the largest power of $p$ that divides $x$, and note that $b_1(x)p=(p+1)^x-1$. Multiplication by $p+1$ can be viewed as an automorphism of $\ZZ_{p^{f+1}}$, and it was shown in~\cite{KovacsNedela2017} that the order of this automorphism is $p^f$. It follows that $(p+1)^x \equiv 1 \pmod{p^{f+1}}$, and consequently $b_1(x)p\equiv (p+1)^x-1 \equiv 0 \pmod{p^{f+1}}$. This proves that $p^f$ divides $b_1(x)$. On the other hand, since the order of $b_1$ is finite, we deduce that $p^{f+1}$ cannot divide $b_1(x)$, for this would force $p^{f+1}$ to divide $x$, which is impossible. This proves that the largest power of $p$ that divides $b_1(x)$ is indeed $p^f$. We will also use the fact (which is a straightforward consequence of the definition of $b_1$) that $b_1(x+y)=b_1(x)+(p+1)^xb_1(y)$ for each $x,y\in \ZZ_{p^e}$.  

As the next step, we show that if $e\geq 2$, then $s_{1,1,1,p^{e-2}}$ has kernel of order $p$. To simplify notation let $s_e=s_{1,1,1,p^{e-2}}$. First consider the case $e=2$. Then $s_2=s_{1,1,1,1}$, and we claim that this is a proper skew morphism of $\ZZ_{p^2}$. In~\cite{KovacsNedela2017} it was proved that two distinct $4$-tuples (that satisfy conditions \ref{cond:1}, \ref{cond:2} and \ref{cond:3}) give two distinct skew morphisms. It can be easily seen that there are exactly $p(p-1)$ skew morphisms of the form $s_{i,0,k,0}$ and that they are all automorphisms of $\ZZ_{p^2}$. Since $|\Aut(\ZZ_{p^2})|=p(p-1)$, it follows that all other skew morphisms of $\ZZ_{p^2}$ must be proper. Hence $s_2$ is proper, and so its kernel must be a non-trivial proper subgroup of $\ZZ_{p^2}$, which implies $|\ker s_2|=p$. From now on we assume that $e\geq 3$. In particular, we will repeatedly use the fact that $p^{e-2}$ is a multiple of $p$. Since $\ker s_e$ is a non-trivial subgroup of a cyclic $p$-group, in order to show that $|\ker s_e|=p$ it is sufficient to show that $|\ker s_e|<p^2$. This is equivalent with $p^{e-2}\notin \ker s_e$. Suppose to the contrary that $p^{e-2}\in \ker s_e$, and so $s_e(1+p^{e-2})=s_e(p^{e-2}+1)=s_e(p^{e-2})+s_e(1)$. In particular, it follows that $b_1(s_e(1+p^{e-2}))=b_1(s_e(p^{e-2})+s_e(1))$. Recall that $b$ is an automorphism of $\ZZ_{p^e}$ of order $p-1$, and so there exists $t\in \ZZ_{p^e}$ such that $\gcd(t,p)=1$, $b(x)=tx$ for all $x\in \ZZ_{p^e}$, and the multiplicative order of $t$ modulo $p^e$ is $p-1$. Using the definition of $s_e$ and the identities proved in the previous paragraph we obtain the following:
\begin{align*}
b_1\left(s_e(1+p^{e-2})\right) &= b_1b_1^{\, -1}abb_{p^{e-2}}b_1\left(1+p^{e-2}\right)=abb_{p^{e-2}}\left(b_1(1)+(p+1)b_1(p^{e-2})\right) \\
&= abb_{p^{e-2}}\left(1+(p+1)b_1(p^{e-2})\right) = ab\left(1+(p+1)b_1(p^{e-2})\right) \\
&=ab(1)+ab\left((p+1)b_1(p^{e-2})\right)=ab(1)+t(p+1)^2b_1(p^{e-2})\, .
\end{align*}   
Note that we used the fact that $b_1(p^{e-2})$ is a multiple of $p$, and so $b_{p^{e-2}}$ fixes the element $1+(p+1)b_1(p^{e-2})$. We also have
\begin{align*}
b_1\left(s_e(p^{e-2})+s_e(1)\right)&=b_1\left(s_e(1)+s_e(p^{e-2})\right)=b_1(s_e(1))+(p+1)^{s_e(1)}b_1(s_e(p^{e-2}))\\
&= abb_{p^{e-2}}b_1(1) + (p+1)^{s_e(1)} abb_{p^{e-2}}b_1(p^{e-2}) = ab(1) + (p+1)^{s_e(1)}abb_{p^{e-2}}(b_1(p^{e-2})) \\
&= ab(1) + (p+1)^{s_e(1)}ab(b_1(p^{e-2})) = ab(1) +t(p+1)^{s_e(1)+1}b_1(p^{e-2})\, .
\end{align*}
It follows that $t(p+1)^2b_1(p^{e-2})$ is equal to $t(p+1)^{s_e(1)+1}b_1(p^{e-2})$ in $\ZZ_{p^e}$. Since the largest power of $p$ that divides $b_1(p^{e-2})$ is $p^{e-2}$, this is equivalent with $t(p+1)^2 \equiv t(p+1)^{s_e(1)+1} \pmod{p^2}$. Moreover, since $t$ is relatively prime to $p$ and the multiplicative order of $p+1$ modulo $p^2$ is $p$, this is true if and only if $2\equiv s_e(1)+1 \pmod{p}$, and so we have $s_e(1) \equiv 1 \pmod{p}$. In particular we have $s_e(1)=1+px$ for some $x\in \ZZ_{p^{e}}$. On the other hand we have $s_e(1)=b_1^{\, -1}abb_{p^{e-2}}b_1(1)=b_1^{\, -1}ab(1)$, and so we have $b_1(1+px)=ab(1)$. The left hand side is equal to $1+(p+1)b_1(px)$, which is congruent to $1$ modulo $p$. The right hand side is equal to $t(p+1)$, which is congruent to $t$ modulo $p$, therefore $t$ must be congruent to $1$ modulo $p$. Recall that the multiplicative order of $p+1$ modulo $p^e$ is equal to $p^{e-1}$. Clearly all powers of $p+1$ are congruent to $1$ modulo $p$, and it follows that every element of $\ZZ_{p^e}$ that is congruent to $1$ modulo $p$ is some power of $p+1$. Consequently the multiplicative order modulo $p^e$ of every such element is a power of $p$. But $t$ has multiplicative order $p-1$, and so it cannot be congruent to $1$ modulo $p$, contradiction. Hence $p^{e-2}\notin \ker s_e$, and it follows that $|\ker s_e|=p$.

We are now ready to show that the complexity of $s_e$ is $2e-1$. We will use induction on $e$. If $e=2$, then $s_2$ is a skew morphism of $\ZZ_{p^2}$ with kernel of order $p$, and so it must be a proper skew morphism. Then by Proposition~\ref{prop:psquare} we deduce that the complexity of $s_2$ is $3$, which is equal to $2e-1$. Next, let $s_e$ be a skew morphism of $p^e$ for some $e\geq 3$. We have $|\ker s_e|=p$, and so $s^*$ (which can be viewed as $s_e$ taken modulo $p^e/|\ker s_e|$) is a skew morphism of $\ZZ_{p^{e-1}}$. We claim that $s^*$ is the skew morphism $s_{e-1}$ of $\ZZ_{p^{e-1}}$. It can be easily seen that the permutations $b_j$ and $a$ (of $\ZZ_{p^e}$) taken modulo $p^{e-1}$ are exactly the permutations $b_j$ and $a$ of $\ZZ_{p^{e-1}}$. To show that $b$ taken modulo $p^{e-1}$ is an automorphism of order $p-1$, recall that $b$ is given by multiplication by $t$, and the multiplicative order of $t$ is $p-1$. We have shown that if an element of $\ZZ_{p^e}$ is congruent to $1$ modulo $p$, then its multiplicative order must be a power of $p$. It follows that the elements $p^i$ for $i\in \{1, \dots, p-1\}$ are all mutually distinct modulo $p$, and so they are also mutually distinct modulo $p^{e-1}$. In particular, $b$ taken modulo $p^{e-1}$ is an automorphism of order $p-1$, which proves our claim. Hence, using the inductive hypothesis we find that the complexity of $s^*$ is $2(e-1)-1$, and hence the complexity of $s_e$ is $2e-1$.
\end{proof}

The following fact about skew morphisms of cyclic groups follows immediately from Theorem~\ref{thm:unbounded}.

\begin{corollary}
For each non-negative integer $c$ there exists infinitely many skew morphisms of finite cyclic groups with complexity $c$. 
\end{corollary} 
\begin{proof}
By Theorem~\ref{thm:unbounded} we find that for each odd prime $p$ and odd $c$ there exists a skew morphism $\varphi_{p,c}$ of $\ZZ_{p^{(c+1)/2}}$ with complexity $c$.  The assertion for all odd values of $c$ now follows from the fact that there are infinitely many odd primes. Next, we prove the assertion for all even values of $c$. The assertion is clearly true for $c=0$, and so we may assume that $c\geq 2$. Note that $\varphi_{p,c+1}$ is a skew morphism of a cyclic $p$-group with complexity $c+1$, and so (by Proposition~\ref{prop:quo:b} and the definition of the derived skew morphism) we deduce that $(\varphi_{p,c+1})'$ is a skew morphism with complexity $c$, and its order is a non-trivial power of $p$. It follows that each $p$ yields a unique skew morphism with complexity $c$, and so there are infinitely many such skew morphisms.       
\end{proof}

Another consequence of Theorem~\ref{thm:unbounded} is the following observation.

\begin{corollary}\label{cor:comp3}
If $\Comp(\ZZ_{n})\subseteq \{0,1,2,3\}$ and $p$ is an odd prime, then $n$ is not divisible by $p^3$. 
\end{corollary}
\begin{proof}
Suppose to the contrary that there exists a positive integer $n$ and an odd prime $p$ such that $\Comp(\ZZ_{n})\subseteq \{0,1,2,3\}$ and $p^3$ divides $n$. Let $e$ be the largest power of $p$ that divides $n$, and note that $\ZZ_{n}\cong \ZZ_{p^e}\times \ZZ_{m}$ where $m$ is a positive integer such that $n=p^em$. We also let $\iota$ denote the trivial permutation of $\ZZ_m$. Since $e\geq 3$, by Theorem~\ref{thm:unbounded} we deduce that there exists a skew morphism $s$ of $\ZZ_{p^e}$ with complexity at least $5$. It can be easily checked that $s \times \iota$ is a skew morphism of $\ZZ_{p^e}\times \ZZ_{m}$ of order $\ord(s)$, and its power function is given by $\pi_{s \times \iota}(a,b)=\pi_s(a)$ for all $(a,b)\in \ZZ_{p^e}\times \ZZ_{m}$; see~\cite[Lemma~7.2]{ConderJajcayTucker2016}. Then $(s \times \iota)'$ is a skew morphism of $\ZZ_{\ord(s)}$ given by $(s \times \iota)'(a)=\sigma_{s \times \iota}(a,(1,1))=\sum_{0\leq i \leq a-1}\pi_{s \times \iota}((s \times \iota)^i(1,1))=\sum_{0\leq i \leq a-1}\pi_{s \times \iota}(s^i(1),1)=\sum_{0\leq i \leq a-1}\pi_s(s^i(1))=\sigma_s(a,1)=s'(a)$ for each $a\in \ZZ_{\ord(s)}$. It follows that $(s \times \iota)'$ is equal to $s'$, and so $s \times \iota$ has the same complexity as $s$. But then the complexity of $s \times \iota$ is at least $5$, which is impossible.  
\end{proof}

Note that Corollary~\ref{cor:comp3} implies that if the order of a cyclic group is divisible by a third power of an odd prime, then it admits a skew morphism of complexity at least $4$. The opposite, however, is not true, even if we restrict ourselves to cyclic groups of odd orders. For example, the cyclic group of order $147$ admits a skew morphism with complexity $4$; see~\cite{census}.

\subsection{Skew morphisms of given order}%%%%%%    

As we have seen earlier, skew morphisms of $\ZZ_n$ have been classified for various infinite families of values of $n$. For example, skew morphisms of $\ZZ_n$ are fully classified if $n$ and $\phi(n)$ are relatively prime, if $n=2^em$ with $e\in\{0,1,2,3,4\}$ and $m$ odd and square-free, or if $n$ is a power of an odd prime. In this section we will not restrict ourselves to specific orders of underlying groups, instead we will restrict the order of skew morphisms. The following theorem characterises and enumerates all proper skew morphism of prime orders.

\begin{theorem}\label{thm:prime}
If $p$ is a prime, then $\ZZ_n$ admits a proper skew morphism of order $p$ if and only if $p$ divides $n$, and $\gcd(p-1,n)\geq 2$. The skew morphisms of $\ZZ_n$ of order $p$ are exactly the permutations given by $a \mapsto a+u(1+v+\dots+v^{a-1})$ where $u\in \{\frac{n}{p}, \frac{2n}{p},\dots, \frac{(p-1)n}{p}\}$ and $v$ is an element from $\{2, 3,\dots,p-1\}$ such that its multiplicative order modulo $p$ divides $n$. The number of such skew morphisms is given by the following formula:
\[ \left(\sum_{\substack{k \mid p-1,\, k \mid n \\ k\neq 1}} \phi(k)\right) (p-1) \, . \]
In particular, $\ZZ_n$ has at most $p^2-3p+2$ proper skew morphisms of order $p$ and the equality holds if and only if $(p-1)p$ divides $n$. 
\end{theorem}
\begin{proof}
By Proposition~\ref{prop:quo:f} we know that there are no proper skew morphisms of $\ZZ_n$ of order $2$, and we also have $\gcd(p-1,n)=\gcd(1,n)=1$, and so the assertion is clearly true for $p=2$. Next, assume that $p\geq 3$. Let $\varphi$ be a proper skew morphism of $\ZZ_n$ such that $\ord(\varphi)=p$. Then by Theorem~\ref{thm:orderofskewcyclic} we deduce that $p$ divides $n$. Next, note that $\varphi'$ is a skew morphism of $\ZZ_p$, and since $\ZZ_p$ admits no proper skew morphisms, we deduce that the complexity of $\varphi'$ is at most $1$. Consequently the complexity $c$ of $\varphi$ is at most $2$, and since $\varphi$ is a proper skew morphism, we have $c=2$. It follows that the auto-order $m$ of $\varphi$ is equal to the order of $\varphi'$, and so by Proposition~\ref{prop:quo:b} we have $m=n/|\ker\varphi|$. This implies that $m$ divides $n$ and $\ker\varphi=\langle m\rangle$. Moreover, since $m$ is equal to the order of an automorphism of $\ZZ_p$, we deduce that $m$ divides $p-1$. We will show that if $m$ is not equal to $1$ and divides both $n$ and $p-1$, then there are exactly $\phi(m)(p-1)$ possible choices for $\varphi(1)$, $\varphi'$ and $\varphi\vert_{\langle m\rangle}$ such that this triple satisfy all five properties listed in Proposition~\ref{prop:even}. Since the complexity of $\varphi$ is even, this will count all skew morphisms of $\ZZ_n$ satisfying $\ord(\varphi)=p$ and $\ord(\varphi')=m$.

It is well known that the number of automorphisms of $\ZZ_p$ of order $m$ (where $m$ divides $p-1$) is exactly $\phi(m)$, and so we have exactly $\phi(m)$ possible candidates for $\varphi'$. Next, since $m\in \langle m\rangle = \ker\varphi$, we have $\varphi(m)+\varphi(1)=\varphi(m+1)=\varphi(1)+\varphi^{\pi_{\varphi}(1)}(m)$, and consequently $\varphi(m)=\varphi^{\pi_{\varphi}(1)}(m)$. This can be rewritten as $\varphi^{\pi_{\varphi}(1)-1}(m)=m$. We know that $\varphi$ is a proper skew morphism of order $p$, and so $\pi_{\varphi}(1)-1$ is an integer between $1$ and $p-2$. Moreover, (again by the fact that $\ord(\varphi)=p$) we have $\varphi^p(m)=m$, and consequently $\varphi^p(m)=\varphi^{\pi_{\varphi}(1)-1}(m)$. Since $p$ and $\pi_{\varphi}(1)-1$ are clearly relatively prime, we deduce that $m$ is fixed by $\varphi$. Hence by Theorem~\ref{thm:orderoforbit} it follows that $\varphi\vert_{\langle m\rangle}$ is the trivial permutation. The last step is to identify all candidates for $\varphi(1)$. Since the complexity of $\varphi$ is even, we know that $(\varphi(1)-1)/m$ is equal to some integer $f$, and so $\varphi(1)=1+fm$. Then, since $\varphi\vert_{\langle m\rangle}$ is trivial, by Equation~\eqref{eq:even} we have $\varphi^e(1)=1+efm$ for each positive integer $e$. But then $1=\varphi^p(1)=1+pfm$, and so $fm$ must be a multiple of $n/p$. There are $p$ distinct multiples of $(n/p)$ in $\ZZ_n$, but we exclude the zero multiple since this results in the trivial skew morphism of $\ZZ_n$. In total there are $p-1$ choices for $\varphi(1)$, $\phi(m)$ choices for $\varphi'$ and a single choice for $\varphi\vert_{\langle m\rangle}$. We will now show that every such triple gives a skew morphism of $\ZZ_n$.

Let $a$ be any integer in $\{1,2,\dots,p-1\}$, and let $h=1+an/p$. Further, let $t$ be any integer such that multiplication by $t$ defines an automorphism $\alpha$ of $\ZZ_p$ of order $m$ (so, in particular, we have $t^m\equiv 1 \pmod{p}$), and let $\beta$ be the trivial permutation of $\ZZ_{n/m}$. Also assume that $m$ is at least two and divides both $p-1$ and $n$, and that $p$ divides $n$. We will show that the triple $(h,\alpha,\beta)$ satisfy all five properties listed in Proposition~\ref{prop:even}, and so it defines a skew morphism of $\ZZ_n$. Let $h'=h-1=an/p$, so $\psi^x(1)=1+xh'$ for each positive integer $x$. Since $m$ divides both $n$ and $p-1$, we deduce that $m$ divides $h'$, and it follows that $\psi^x(1)\equiv 1 \pmod{m}$. Since $\alpha$ is an automorphism, we also have $\pi_{\alpha}(x)=1$, and so \ref{prop:even:a} holds. Next, \ref{prop:even:b} is true as both $\beta'$ and $\alpha^m$ are trivial permutations. We proceed to verify \ref{prop:even:c}. Since $\beta$ is trivial, the right hand side is equal to $xm$. From Equation~\eqref{eq:even2} we know that the left hand side is equal to $(1+h')+(1+\alpha(1)h')+(1+\alpha^2(1)h')+\dots+(1+\alpha^{xm-1}(1)h')$. Since $\alpha$ is an automorphism of order $m$, this can be rewritten as $xm+h'(1+\alpha(1)+\dots+\alpha^{m-1}(1))x$. Also note that $(1+\alpha(1)+\dots+\alpha^{m-1}(1))(t-1) \equiv (1+t+\dots+t^{m-1})(t-1) \equiv t^m-1 \equiv 0 \pmod{p}$, and so $p$ divides $1+\alpha(1)+\dots+\alpha^{m-1}(1)$. It follows that $n$ divides $h'(1+\alpha(1)+\dots+\alpha^{m-1}(1))$, therefore $\varphi(xm)=xm$. Next, using some of the above observations, \ref{prop:even:d} can be easily verified as follows: 
\begin{align*}
\varphi(x+m) &= (1+h')+(1+th')+(1+t^2h')+\dots+(1+t^{x+m-1}h') \\ 
&=(1+h')+\dots+(1+t^{x-1}h')+(1+t^{x}h')+\dots+(1+t^{x+m-1}h') \\
&=\varphi(x)+(1+t^{x}h')+\dots+(1+t^{x+m-1}h') \\
&=\varphi(x)+m+h't^x(1+t+\dots+t^{m-1}) \\
&=\varphi(x)+m+0=\varphi(x)+\beta^{\alpha^x(1)}(1)m=\varphi(x)+\varphi^{\alpha^x(1)}(m)\, .
\end{align*} 
Finally, since $\psi^x(1)=1+xh'$ with $h'=an/p$, and $a$ is not divisible by $p$, we deduce that the smallest positive integer $i$ such that $\psi^i(1)=1$ is $p$, and so \ref{prop:even:e} holds.

Since $p$ divides $n$ we have $n\geq p$, and so by Lemma~\ref{lem:comp:a} it follows that $m\geq 2$. For all such $m$ we have shown that $\varphi$ exists if and only if $m$ divides both $p-1$ and $n$, and the number of corresponding skew morphisms is $\phi(m)(p-1)$. This proves the second assertion. To prove the third assertion note that $\sum_{k \mid p-1} \phi(k) = p-1$. It follows that    
\[ \sum_{\substack{k \mid p-1,\, k \mid n \\ k\neq 1}} \phi(k)(p-1) \leq \sum_{\substack{k \mid p-1 \\ k\neq 1}} \phi(k)(p-1) = \left(\sum_{k \mid p-1} \phi(k)(p-1)\right)-\phi(1)(p-1)= (p-1)^2-(p-1) = p^2-3p+2\, ,   \]
and the equation holds if and only if all divisors of $p-1$ also divide $n$. This is clearly true if and only if $p-1$ divides $n$, and the rest follows easily. 

To show that every skew morphism of $\ZZ_n$ has the form given in the statement of the theorem, recall that $\varphi(a)=(1+h')+(1+\alpha(1)h')+(1+\alpha^2(1)h')+\dots+(1+\alpha^{a-1}(1)h')=a+h'(1+\alpha(1)+\dots+\alpha^{a-1}(1))=a+h'(1+t+\dots+t^{a-1}(1))$. The claim now follows easily from the fact that $h'$ is a non-zero multiple of $(n/p)$ in $\ZZ_n$ and $t$ is an unit modulo $p$ with multiplicative order $m$. 
\end{proof}

An interesting consequence of Theorem~\ref{thm:prime} is that for each prime $p$ there exists a constant $N_p$ such that for every positive integer $n$ the number of proper skew morphisms of $\ZZ_n$ is bounded above by $N_p$. (In particular, we have $N_p=p^2-3p+2$.) This is actually quite surprising, since as we prove next this is not true for automorphisms of cyclic groups.

\begin{proposition}\label{prop:unbounded}
For each pair of integers $\ell$ and $N$ such that $\ell\geq 2$ there exists a finite cyclic group that admits at least $N$ automorphisms of order $\ell$. 
\end{proposition}
\begin{proof}
By Dirichlet’s theorem there exist $N$ distinct primes $p_1, p_2, \dots, p_N$ such that each of them is congruent to $1$ modulo $\ell$. Let $n=p_1p_2\dots p_N$ and let $G$ be the cyclic group of order $n$. It is well known that $\Aut(G)$ is isomorphic to the multiplicative group of integers modulo $n$, and hence in our case we have $\Aut(G) \cong \ZZ_{p_1-1} \times \ZZ_{p_2-1} \times \dots \times \ZZ_{p_N-1}$. For each $i\in \{1,2,\dots,N\}$ the order of the cyclic factor $\ZZ_{p_i-1}$ is divisible by $\ell$, and so it contains an element $a_i$ of order $\ell$. Now for every $i$ we define an element $b_i$ of $\Aut(G)$ as follows:
\[    (b_i)_j= 
\begin{cases}
    a_i & \text{if } j=i,\\
    0,              & \text{otherwise.}
\end{cases}\]
It can be easily seen that elements $b_i$ are mutually distinct automorphisms of $G$ of order $\ell$, which proves our claim.  
\end{proof}

Proposition~\ref{prop:unbounded} implies that for any $\ell\geq 2$ the number of automorphisms (of order $\ell$) of a finite cyclic group can be arbitrarily large. In contrast, the number of proper skew morphisms (of order $\ell$) of a finite cyclic group is bounded for some values of $\ell$. If this is the case we say that the integer $\ell$ is \emph{skew bounded}. Clearly $\ell=1$ is skew bounded and by Theorem~\ref{thm:prime} all primes are skew bounded. An interesting question is if there are any other skew bounded positive integers. Here we look at the case $\ell=4$. Similar to the prime case we know that all skew morphisms of $\ZZ_4$ are group automorphisms. Moreover, if $\varphi$ is a skew morphism of order $4$, then there are only two possible candidates for its reduction $\varphi'$. Somewhat surprisingly, it turns out that $4$ is not skew bounded. The key difference is that while in the prime case the reduction $\varphi\vert_{\langle m\rangle}$ must be trivial, in the case when $\ord(\varphi)=4$ the reduction $\varphi\vert_{\langle m\rangle}$ can be non-trivial.

\begin{theorem}
Let $n=16p_1p_2\dots p_i$ where $p_1$, $p_2$, $\dots$, $p_i$ are mutually distinct odd primes. Then the number of proper skew morphisms of $\ZZ_n$ of order $4$ is $2^{i+2}$. In particular, $4$ is not skew bounded.  
\end{theorem}
\begin{proof}
Let $\varphi$ be a proper skew morphism of $\ZZ_n$ of order $4$ with complexity $c$ and auto-order $m$. Then by Proposition~\ref{prop:quo:f} we know that $\varphi'$ must be a non-trivial skew morphism of $\ZZ_4$, and so according to~\cite{ConderList} we find that $\varphi'$ is the automorphism $(1,3)$ of $\ZZ_4$. Then, since $2=\ord(\varphi')=n/|\ker\varphi|$, we deduce that $\ker\varphi$ is the subgroup of $\ZZ_n$ generated by the element $2$. Furthermore, by Lemma~\ref{lem:comp:f} we have $c=2$, and by Proposition~\ref{prop:quo:b} we have $\pi_{\varphi}(1)=3$. Since $\ker\varphi = \langle 2\rangle$, it follows that $\varphi\vert_{\langle 2\rangle}$ is an automorphism of $\ZZ_{n/2}$. Let $s$ be the unit modulo $n/2$ such that $\varphi\vert_{\langle 2\rangle}$ is given by $x\mapsto sx$, and let $h=\varphi(1)-1$. In~\cite{BachratyJajcay2017} it was shown that every skew morphism of $\ZZ_n$ with complexity $2$ is uniquely determined by the set of four parameters: the smallest non-zero element $d'$ of $\ker\varphi$, the element $h'=\varphi(1)-1$, the smallest positive integer $s'$ such that $\varphi(d)=s'd$, and the positive integer $e'$ equal to $\pi_{\varphi}(1)$. In our case we have $d'=2$, $h'=h$, $s'=s$ and $e'=3$. It was also proved that a quadruple of parameters defines a skew morphism if and only they satisfy seven numerical properties. In our case some of them are trivially true, and they can be reduced to the following four properties:  
\begin{enumerate}[label={\rm (P\arabic*)}]
\item\label{item:1} $h$ is a multiple of $2$ strictly smaller than $n$;
\item\label{item:2} $r=4$ is the smallest positive integer such that $h\sum_{i=0}^{r-1}  s^i \equiv 0 \pmod{n}$;
\item\label{item:3} $2s\equiv 2+h(2+s+s^2) \pmod{n}$; and
\item\label{item:4} $s^2 \equiv 1 \pmod{n/2}$.
\end{enumerate}
Thus to count all proper skew morphisms of $\ZZ_n$ of order $4$ we need to count the number of pairs $s$ and $h$ for which properties~\ref{item:1} to~\ref{item:4} are all true. From~\ref{item:1} we know that $h=2f$ for some $f\in \ZZ_{n/2}$, and instead of counting suitable pairs $(s,h)$ we will be finding pairs $(s,f)$. Then we can rewrite~\ref{item:3} as $s\equiv 1+f(2+s+s^2) \pmod{n/2}$, and using~\ref{item:4} this can be further rewritten as 
\begin{equation}\label{eq:four}
s\equiv 1+f(3+s) \pmod{n/2}\,.
\end{equation}
Next, from~\ref{item:4} we know that $s$ is a unit modulo $n/2$ of order at most $2$. Since $n/2=8p_1p_2\dots p_i$, it follows that $s$ is also a unit of order at most $2$ modulo $8$, and also a unit of order at most $2$ modulo $p_j$ for each $j\in \{1,2,\dots,i\}$. The only units of order at most $2$ modulo a prime are $-1$ and $1$, so for each $j$ we have two possibilities. If $s \equiv 1 \pmod{p_j}$, then Equation~\eqref{eq:four} yields $0\equiv 4f \pmod{p_j}$, and so $f \equiv 0 \pmod{p_j}$. On the other hand, if $s \equiv -1 \pmod{p_j}$, then Equation~\eqref{eq:four} yields $-2\equiv 2f \pmod{p_j}$, and so $f \equiv -1 \pmod{p_j}$. Note that~\ref{item:2} implies that $f(1+s+s^2+s^3)\equiv 0 \pmod{p_j}$, and by plugging in~\ref{item:4} this yields $2f(1+s)\equiv 0 \pmod{p_j}$. This is always true, since $p_j$ always divides either $f$ or $1+s$. Note that the same argument can be also applied to show that $f(1+s)\equiv 0 \pmod{p_j}$, and so since $r\neq 2$ we must have $f(1+s)\not\equiv 0 \pmod{8}$. The units of order at most $2$ modulo $8$ are exactly $1$, $3$, $5$ and $7$, and the respective equations derived from Equation~\ref{eq:four} are $0\equiv 4f$, $2\equiv 6f$, $4\equiv 0$ and $6\equiv 2f$, all taken modulo $8$. This gives $8$ possible candidates for the pair $(s,f) \bmod{8}$, namely $(1,0)$, $(1,2)$, $(1,4)$, $(1,6)$, $(3,3)$, $(3,7)$, $(7,3)$ and $(7,7)$. We argued earlier that $f(1+s)$ cannot be congruent to $0$ modulo $8$, and it can be easily checked that this will exclude the pairs $(1,0)$, $(1,4)$, $(7,3)$ and $(7,7)$. This reduces suitable pairs to the following four: $(1,2)$, $(1,6)$, $(3,3)$ and $(3,7)$. For each of them it can be easily checked that the smallest $r$ in~\ref{item:4} is indeed $4$. The above analysis and Chinese Remainder Theorem now implies that we have four choices for $(s,f) \bmod{8}$ and two choices for $(s,f) \bmod{p_j}$, any specific choice gives a unique pair $(s,f) \bmod{n/2}$ such that $s$ and $h=2f$ satisfy~\ref{item:1}, \ref{item:2}, \ref{item:3} and \ref{item:4}, and there are no other suitable pairs $(s,f)$. It follows, that the number of proper skew morphisms of $\ZZ_n$ is exactly $2^{i+2}$.
\end{proof}

In summary, we have shown that all non-composite positive integers are skew bounded, while $4$ is not skew bounded. We leave the problem of determining whether there exists a composite skew-bounded positive integer an open question.  

%%%%%%%%%%%
\section{Remarks}\label{sec:remarks}
%%%%%%%%%%%

In Section~\ref{sec:main} we provided a characterisation of skew morphisms for all finite cyclic groups. As illustrated in Sections~\ref{sec:recs} and~\ref{sec:apps} this can be used to generate a census of skew morphisms of finite cyclic groups up to any given order, and also to prove various theoretical observations. On the other hand, there are certain limitations to this characterisation, mainly due to its recursive nature. Most notably, there is no obvious way how it can be used to find an explicit formula for the number of skew morphisms of $\ZZ_n$. In fact, there is no known closed formula even in a much simpler case when we restrict ourselves to skew morphisms of $\ZZ_n$ of complexity $2$. Interestingly, Theorem~\ref{thm:unbounded} suggests that there are skew morphisms of arbitrarily large complexity that have also an alternative simpler description. In particular, as proved in~\cite{KovacsNedela2017}, every skew morphism of $\ZZ_{p^e}$ (where $p$ is an odd prime) can be constructed as a composition of at most five permutations of $\ZZ_{p^e}$, even though the set of complexities of all such skew morphisms is unbounded. For the above reasons, we consider the search for an alternative characterisation of skew morphisms of finite cyclic groups to be worth further investigation.

An interesting problem in the context of derived skew morphisms is to determine whether each skew morphism of a finite cyclic group is a derived skew morphism of some other skew morphism of a finite cyclic group. If we consider the set $S_{2000}$ of all skew morphisms of cyclic groups of order up to $2000$ (obtained from~\cite{census}), then for each skew morphism $\rho$ of $\ZZ_n$ with $n\leq 46$ there exists $\varphi\in S_{2000}$ such that $\rho=\varphi'$. In the case when $n=47$ there exists a skew morphism $\ZZ_n$ which cannot be obtained as a derived skew morphism of a skew morphism from $S_{2000}$, namely any automorphism of $\ZZ_{47}$ of order $46$. But it can be easily checked that for each odd prime $p$ an automorphism of $\ZZ_p$ of order $p-1$ is a derived skew morphism of a suitable skew morphism of the cyclic group of order $p(p-1)$. Note that the smallest prime $p$ for which $p(p-1)>2000$ is indeed $47$. These observations suggest that if we take the set $S$ of all skew morphisms of finite cyclic groups, than it might be true that for each $\rho\in S$ there exists $\varphi\in S$ such that $\rho=\varphi'$. We leave this as an open question. One might even ask if for each $\rho\in S$ there are infinitely many $\varphi\in S$ such that $\rho=\varphi'$.

\providecommand{\bysame}{\leavevmode\hbox to3em{\hrulefill}\thinspace}
\providecommand{\MR}{\relax\ifhmode\unskip\space\fi MR }
\providecommand{\MRhref}[2]{%
  \href{http://www.ams.org/mathscinet-getitem?mr=#1}{#2}
}
\providecommand{\href}[2]{#2}

\end{document}